\documentclass[12pt]{amsart}
\usepackage{bm}
\usepackage{amssymb,amsthm,amsmath,amscd}
\usepackage[numbers,sort&compress]{natbib}
\usepackage{color}
\usepackage{graphicx}
\usepackage{epsf,epsfig}
\usepackage{epsfig,multicol}
\usepackage{tikz,tikz-cd}
\usepackage{hyperref}
\usepackage{enumitem}
\usepackage{booktabs}
\usepackage[all]{xy}

\usepackage{enumitem}
\setlist[itemize]{label=\textbullet}

 \newcommand{\Z}{\mathbb{Z}}
 
 \newcommand{\Zc}{\mathcal{Z}}

 \newcommand{\Ac}{\mathcal{A}}

\def\im{\mathop{\rm Im}\nolimits}

\def\Hom{{\mathrm{Hom }}}

\newtheorem{problem}{Problem}

\usepackage{geometry}

\title{Representation characterization of  equivariant geometric bordism}
\author{Hao Li, Bo Chen, Zhi L\"{u} and  Qifan Shen}

\address{College of Mathematics and Sciences,  Shanghai Normal University,  Shanghai, 200234, 
P. R. China}
\email{14110840001@fudan.edu.cn}

\address{School of Mathematics and Statistics, Huazhong University of Science and Technology, Wuhan, China} 

\email{bobchen@hust.edu.cn}

\address{School of Mathematical Sciences, Fudan University, Shanghai,  200433, P. R. China}
\email{zlu@fudan.edu.cn}

\address{Shanghai High School, Shanghai, 200231, P. R. China}

\email{17110180009@fudan.edu.cn}

\thanks{ Partially supported by the grant from NSFC11971112, NSFC12501087}

\theoremstyle{plain}

\newtheorem{theorem}{Theorem}[section]

\newtheorem{lemma}[theorem]{Lemma}

\theoremstyle{definition}
\newtheorem{definition}[theorem]{Definition}

\newtheorem{example}{Example}
\newtheorem{remark}{Remark}

\allowdisplaybreaks
\begin{document}

\begin{abstract}
Let $G_k=(\mathbb Z_2)^k$.  We establish a representation characterization determining when finitely many faithful \(G_k\)-representations can serve as
the fixed-point data of a smooth closed
$G_k$-manifold. By partitioning representations via multiplicities of nonzero weights \(\rho\) and isomorphic restrictions to \(\ker\rho\), we formulate an evenness condition; realizability is equivalent to this condition, yielding a finite local restatement of the tom Dieck–Kosniowski–Stong integrality criterion.
As an application, we prove that \(\dim_{\mathbb Z_2}\mathcal Z_4(G_3)=32\) and exhibit 32 basis elements represented by \(G_3\)-actions on \(\mathbb RP^2\times\mathbb RP^2\), \(\mathbb RP^4\) and \(\mathbb RP(\gamma \oplus \gamma \oplus \gamma  \oplus\underline{\mathbb R})\), obtained from three basic actions by precomposition with automorphisms of \(G_3\).
\end{abstract}
\maketitle

\section{Introduction}

In the 1960s, Conner and Floyd \cite{CF} initiated the study of
equivariant geometric bordism while investigating differentiable periodic maps
on smooth manifolds.  In 1970, tom Dieck \cite{tom1} introduced
homotopy-theoretic equivariant (co)bordism.  Since then, various
geometric and homotopy-theoretic equivariant (co)bordism theories have
been widely studied.

\vskip.2cm
Throughout this paper we write $G_k=(\Z_2)^k$.  A typical geometric
equivariant bordism ring is $\mathfrak{N}^{G_k}_*$, which consists of
equivariant bordism classes of smooth closed $G_k$-manifolds, with no
effectivity assumption of $G_k$-actions.  It  has a natural  module structure over the ordinary
unoriented bordism ring $\mathfrak{N}_*$.

When $k=1$, Conner and Floyd \cite{CF} proved that the equivariant
bordism class of a smooth closed $G_1$-manifold is determined by its
fixed-point components together with their equivariant normal bundles.
Stong \cite{St} and tom Dieck \cite{d} extended this fixed-point
detection principle to $G_k$-actions.

Beyond $k=1$, the structure of $\mathfrak{N}^{G_k}_*$ remains far from
fully understood; see \cite{Al,Sinha} for $\mathfrak{N}^{G_1}_*$.
Conner and Floyd \cite{CF}
focused attention on the subring $\mathcal{Z}_*(G_k)$ formed by classes of smooth
closed $G_k$-manifolds with isolated fixed points.  They introduced
the \emph{Conner--Floyd unoriented $G_k$-representation algebra}
$\mathcal{R}_*(G_k)$ and defined the homomorphism
\[
 \phi_*:\ \mathcal{Z}_*(G_k)\longrightarrow\mathcal{R}_*(G_k)
\]
sending an equivariant bordism class $[M]$ to the sum
$\sum_{p\in M^{G_k}}\tau_pM$ of the tangent representations at
its isolated fixed points.  Here $\mathcal{R}_*(G_k)$ is the polynomial
algebra over $\Z_2$ generated by the  irreducible real
representations of $G_k$, in which addition is the formal sum and multiplication is the Whitney sum.
Stong \cite{St} subsequently proved that $\phi_*$ is a monomorphism.  When
$k=1,2$, we have known from Conner--Floyd's work \cite[Theorems 25.1,
31.1 and 31.2]{CF} that $\mathcal{Z}_*(G_1)\cong \mathbb{Z}_2$ and
$\mathcal{Z}_*(G_2)\cong \mathbb{Z}_2[u]$, where $u$ is
$\mathbb{R}P^2$ with standard $G_2$-action.  For $k\geq 3$, the indecomposable elements of $\mathcal{Z}_*(G_k)$ were studied in \cite{M} and \cite{MS}; however,  
the
ring structure of $\mathcal{Z}_*(G_k)$ remains open.

In contrast, the  subring $\mathcal{Z}_*(G_k)$ is not naturally an
$\mathfrak{N}_*$-module: multiplying by a positive-dimensional
manifold with trivial $G_k$-action replaces an isolated fixed point by
a positive-dimensional fixed component.  Effectivity of $G_k$-actions also requires
some care in this setting.  Let $K$ be the kernel of a $G_k$-action on
$M$.  The action factors uniquely through $G_k/K$, and the induced
$G_k/K$-action is effective.  Since the original action and the quotient
action have the same image in $\operatorname{Diff}(M)$, they have the
same fixed-point set.  Moreover, at every fixed point $p$, the tangent
representation of $G_k$ is the pullback, along the quotient map
$q\colon G_k\to G_k/K$, of the tangent representation of $G_k/K$.
Thus, if $K\neq\{1\}$, the entire fixed-point data is obtained from an
effective action of the lower-rank group
$G_k/K\cong G_{k-\dim_{\Z_2}K}$.  Hence non-effective actions introduce
no genuinely new fixed-point data at rank $k$.

For disconnected manifolds there is a further distinction: an action
may be effective on the whole manifold without being effective on
each connected component.  Guan and L\"u \cite{GLN} recently emphasized
this distinction by working with the \emph{fully effective} subring,
in which the action on every connected component is effective.  

Throughout this paper we convention that $\mathcal Z_*(G_k)$ is
understood to be represented by connected manifolds with effective
$G_k$-actions and isolated fixed points.  This  is because  the
disjoint union of two connected $G_k$-manifolds is equivariantly bordant to a connected
$G_k$-manifold obtained by an equivariant connected sum along free orbits.

For a connected manifold with an effective $G_k$-action, the tangent representation
at every fixed point is faithful.  If the fixed point is isolated, it
also has no trivial summand. 
 Since $\phi_*$ is injective,
$\mathcal Z_*(G_k)$ is identified with $\im\phi_*\subseteq
\mathcal R_*(G_k)$.  Determining 
which collections of faithful representations with no trivial summand
are realizable as fixed-point data is therefore a fundamental
problem.



\begin{problem}\label{problem}
Which elements of $\mathcal R_n(G_k)$ occur as the fixed-point data
of smooth closed $n$-dimensional $G_k$-manifolds?
\end{problem}

Two complementary approaches to Problem~\ref{problem} have shaped its
development.  The first is localization in equivariant bordism.  tom
Dieck \cite{d} proved that equivariant characteristic numbers determine
the geometric bordism class and, for actions with isolated fixed
points, characterized geometric fixed-point data by an integrality
condition.  Kosniowski and Stong \cite{KS} subsequently expressed the
relevant characteristic numbers in terms of fixed-point data by a
localization formula.  Together, these results give the following form
of the tom Dieck--Kosniowski--Stong criterion.

\begin{theorem}[tom Dieck; Kosniowski-Stong]\label{Kos-Stong}
Let $\Ac$ be a finite set of faithful $n$-dimensional
$G_k$-representations with no trivial summand.  Then
$\sum_{\tau\in \Ac} \tau \in\im\phi_n$ (equivalently,
$\Ac$ is the fixed-point data of some
$G_k$-manifold) if and only if, for every symmetric polynomial
$f(x_1,\cdots,x_n)$ with coefficients in $\Z_2$,
\[
 \sum_{\tau\in \Ac} \frac{f(\tau)}{\chi^{G_k}(\tau)}
 \in H^*(BG_k;\Z_2)\cong\Z_2[\rho_1,\cdots,\rho_k],
\]
where $\chi^{G_k}(\tau)$ is the product of the $n$ 
weights of $\tau$, and $f(\tau)$ denotes the substitution of these
weights into $f$.
\end{theorem}

The second approach is combinatorial.  Goresky, Kottwitz, and MacPherson
\cite{GKM} related the equivariant cohomology rings of certain
algebraic varieties to regular labeled graphs; see also \cite{GZ}.
This viewpoint has also been applied to $G_k$-manifolds
\cite{BGH,L1,L3} and provides the graph-theoretic language used in the
square-free case.

These approaches have led to several computations and structural
results.  Using the Davis--Januszkiewicz theory of small covers
\cite{DJ}, L\"u \cite{L2} computed $\mathcal{Z}_3(G_3)$; L\"u and Tan
\cite{LT} used a differential on the dual representation algebra to
characterize $\im\phi_k$ for $G_k$: a faithful polynomial is realizable
as fixed-point data if and only if its dual vanishes by the
differential; and Chen, L\"u, and Tan \cite{CLT} gave a homological
description of $\mathcal{Z}_k(G_k)$ and determined its dimension.  In
the square-free case, L\"u
\cite{L1} gave a labeled graph characterization of the fixed-point
data.  Once weights at fixed points are allowed to repeat, however, the
labeled graph description does not by itself capture the additional
multiplicity data.

The present paper solves Problem~\ref{problem} by translating the
integrality criterion in Theorem~\ref{Kos-Stong} into a finite parity
test on the weights of the given representations.  The geometric input
is as follows.  Fix a nonzero weight $\rho$ and
put $H=\ker\rho$.  If $C$ is a component of $M^H$ containing a
$G_k$-fixed point, then, for every $p\in C\cap M^{G_k}$, the tangent
representation $T_pC$ is a direct sum of copies of $\rho$, while normal
weights that differ by $\rho$ have the same restriction to $H$.
This geometric decomposition leads on the representation side to the
mod $\rho$ equivalence classes.

The evenness condition records the cancellations within each such
class: certain lower-degree subproducts must occur with even total
multiplicity.  When the multiplicity of $\rho$ is one, this reduces to
the requirement that each class itself has even cardinality.  In higher
multiplicity, every relevant tuple of degree below that multiplicity
must occur an even total number of times across the class.
The precise finite list of congruences is given in
Definition~\ref{def:tuple-condition}.

Our main result is as follows.

\begin{theorem}\label{A}
Let $\Ac$ be a finite set of faithful $n$-dimensional
$G_k$-representations with no trivial summand.  Then $\Ac$ is
the fixed-point data of a smooth closed $n$-dimensional $G_k$-manifold
if and only if $\Ac$ satisfies the evenness condition.
\end{theorem}

Theorem~\ref{A} is equivalent to the tom
Dieck--Kosniowski--Stong criterion, but it organizes the same
information in a different way.  In place of testing the integrality
of localization sums for all symmetric polynomials, it characterizes
realizability by finitely many local parity conditions on mod $\rho$
equivalence classes.  In the square-free case, these conditions reduce
to the labeled-graph characterization of L\"u \cite{L1}: for every
nonzero weight $\rho$, each mod $\rho$ equivalence class has even
cardinality.

The two implications proceed differently.  For necessity, apply the
Kosniowski--Stong bundle formula \cite{KS} to the equivariant normal
bundle of each component of $M^{\ker\rho}$ described above.  This gives
local integrality relations.  After the components belonging to the
same mod $\rho$ equivalence class are combined,
Lemma~\ref{lem:normal-integrality-implies-evenness} converts these
relations into the evenness condition for that class.  Doing this for
all classes and all $\rho$ gives the global evenness condition.

Conversely, the sufficiency argument is algebraic.  The evenness
condition on each mod $\rho$ equivalence class implies its local
integrality condition by
Lemma~\ref{lem:evenness-implies-local-integrality}.  Summing over the
classes recovers the global integrality condition, and
Theorem~\ref{Kos-Stong} then gives the required geometric realization.

As an application, we compute
\[
 \dim_{\mathbb Z_2}\mathcal Z_4(G_3)=32.
\]
The calculation is organized by the three types of faithful
four-dimensional representations, and the resulting generators are
realized by explicit actions on $\mathbb RP^2\times\mathbb RP^2$,
$\mathbb RP^4$, and
$\mathbb RP(\gamma\oplus \gamma\oplus \gamma\oplus\underline{\mathbb R})$.

An earlier version of this work  appeared as the arXiv
preprint \cite{LLS}.  Guan and L\"u subsequently used the
characterization proposed there in their calculation of
$\mathcal Z_5(G_3)$ \cite{GL}.  The present paper gives a revised proof
of the characterization.

The paper is organized as follows.  Section~\ref{sec:prelim}
introduces the evenness condition.  Section~\ref{proof} proves the
necessity of the condition, and Section~\ref{sec:sufficiency} proves
its sufficiency.  Finally, Section~\ref{sec:z4-application} applies
the criterion to $\mathcal Z_4(G_3)$ and constructs geometric
representatives for the generators.

\section{The evenness condition}\label{sec:prelim}

\subsection{Representations of $G_k$ and the polynomial ring}

All irreducible real representations of $G_k$ are one-dimensional and
are indexed by $\Hom(G_k,\mathbb Z_2)$.  Explicitly, a homomorphism
$\rho\colon G_k\to\mathbb Z_2$ determines the representation on
$\mathbb R$ given by
\[
 g\cdot x=(-1)^{\rho(g)}x.
\]
We call $\rho$ the \emph{weight} of this representation.  Thus nonzero
weights index the nontrivial irreducible real representations of
$G_k$.  The Conner--Floyd representation algebra
$\mathcal R_*(G_k)$ is the graded polynomial algebra over $\mathbb Z_2$
generated by the irreducible representations of $G_k$.

The equivariant Euler class gives a homomorphism
\[
 \chi^{G_k}\colon \mathcal R_*(G_k)\longrightarrow
 H^*(BG_k;\mathbb Z_2).
\]
Choose a basis $\{\rho_1,\dots,\rho_k\}$ of
$\Hom(G_k,\mathbb Z_2)$.  Then
\[
 H^*(BG_k;\mathbb Z_2)=R:=\mathbb Z_2[\rho_1,\dots,\rho_k],
\]
and $\chi^{G_k}$ sends an irreducible representation to the
corresponding linear form in $R$.  If an $n$-dimensional representation
$\tau$ has no trivial summand and has weights
$\rho_1,\ldots,\rho_n$, then
\[
 \chi^{G_k}(\tau)=\rho_1\cdots\rho_n.
\]
Unique factorization in $R$ shows that this monomial determines $\tau$
up to isomorphism.  We therefore encode $\tau$ by its \emph{weight
monomial} $\rho_1\cdots\rho_n$ and use the same symbol for the
representation and its weight monomial.

We use the following notation:
\begin{itemize}
    \item $R_1=\Hom(G_k,\mathbb Z_2)\setminus\{0\}$, the set of
    nonzero weights, viewed also as the degree-one elements of $R$;
    \item $S=\{\prod_{\rho\in R_1}\rho^{a_\rho}\mid
    a_\rho\in\mathbb Z_{\geq0}\}$, including the empty product $1$;
    \item $S_n=\{\tau\in S\mid\deg\tau=n\}$ and
    $S_{<n}=\bigcup_{0\leq i<n}S_i$;
    \item for $\beta\in R_1$,
    $S_\beta=\{\tau\in S\mid\beta\nmid\tau\}$.
\end{itemize}
All fixed point data below are viewed  as elements in the Conner--Floyd
representation algebra with $\mathbb Z_2$ coefficients.  Thus a
finite set $\Ac\subset S_n$ records the terms which occur with
coefficient one after equal tangent representations have been
cancelled in pairs.

Both $S$ and $S_\beta$ are multiplicatively closed subsets of $R$, and $S_\beta^{-1}R \subsetneq S^{-1}R$.
The passage between global and local integrality conditions in the localization theorem relies on the following elementary observation.

\begin{lemma}\label{cap}
    $\bigcap_{\beta \in R_1} S_\beta^{-1}R = R$.
\end{lemma}

\begin{proof}
Every denominator in $S$ is a product of elements of $R_1$.  If
$q\in S^{-1}R$ belongs to $S_\beta^{-1}R$, then the $\beta$-adic
valuation of $q$ is nonnegative.  Thus an element in the displayed
intersection has nonnegative valuation at every possible irreducible
factor of its denominator.  Since $R$ is a unique factorization
domain, it has no denominator and hence lies in $R$.
\end{proof}

\subsection{Multiplicity}
\begin{definition}\label{multiplicity}
For $\tau\in S$ and $\rho\in R_1$, set
\[
 m(\tau,\rho)=\max\{a\in\mathbb Z_{\geq0}\mid\rho^a\mid\tau\}.
\]
A weight monomial $\gamma=\prod_{i=1}^l\beta_i^{a_i}$, where the $\beta_i\in R_1$
are distinct, $a_i>0$, and $\sum_i a_i=t$, is called a
\emph{$t$-tuple}.  Its multiplicity in $\tau$ is
\[
 m(\tau,\gamma)=\prod_{i=1}^l
 \binom{m(\tau,\beta_i)}{a_i}.
\]
Thus $m(\tau,\gamma)$ counts the occurrences of the tuple $\gamma$
among the weights of $\tau$.  We set $m(\tau,1)=1$ for the $0$-tuple.
For a finite set $\Ac\subseteq S$, define
\[
 m(\Ac,\gamma)=\sum_{\tau\in\Ac}m(\tau,\gamma).
\]
When a multiplicity is used as a coefficient of an element of $R$, it
is understood to be reduced modulo $2$.
\end{definition}

\begin{example}\label{ex:occurrence}
    For $k=2$ and $R = \mathbb{Z}_2[\rho_1, \rho_2]$, let $\tau = \rho_1^3 \rho_2^2 \rho_{12}^2 \in S$, where $\rho_{12}=\rho_1+\rho_2$.
    \begin{itemize}
        \item Multiplicities of weights: $m(\tau, \rho_1)=3$, $m(\tau, \rho_2)=2$, $m(\tau, \rho_{12})=2$.
        \item For the $3$-tuple $\gamma_1 = \rho_1^2 \rho_2$,
        $m(\tau, \gamma_1) = \binom{3}{2}\binom{2}{1}  =6. $
        \item For the $3$-tuple $\gamma_2 =  \rho_1\rho_2 \rho_{12}$,
        $ m(\tau, \gamma_2)= \binom{3}{1}\binom{2}{1}\binom{2}{1} =12. $
    \end{itemize}
\end{example}

\subsection{The mod $\rho$ equivalence relation and the evenness condition}
\label{subsec:mod-rho-evenness}

\begin{definition}\label{sim-def}
Fix $\rho\in R_1$.  Weight monomials $\tau,\tau'\in S_n$ are \emph{mod
$\rho$ equivalent}, denoted by $\tau\sim_\rho\tau'$, if
\[
 m(\tau,\rho)=m(\tau',\rho)=d>0
 \quad\text{and}\quad
 \frac{\tau}{\rho^d}\equiv
 \frac{\tau'}{\rho^d}\pmod\rho.
\]
\end{definition}

For later use, suppose that
\[
 \tau=\rho^d\prod_{i=1}^{n-d}\beta_i,
 \qquad
 \tau'=\rho^d\prod_{i=1}^{n-d}\beta'_i,
\]
where $\beta_i,\beta'_i\neq\rho$ for all $i$.  Since $R/(\rho)$ is a
polynomial ring, unique factorization shows that
$\tau\sim_\rho\tau'$ if and only if, after reindexing, each $\beta'_i$
is either $\beta_i$ or $\beta_i+\rho$.  In representation-theoretic
terms, put
$
 H=\ker\rho=\{g\in G_k\mid\rho(g)=0\}.
$
 Two nonzero weights of $G_k$ have the same
restriction to $H$ if and only if they differ by $\rho$.  It follows
that $\tau\sim_\rho\tau'$ if and only if their restrictions to $H$ are
isomorphic.  In particular, $\sim_\rho$ is an equivalence relation on
the weight monomials in $S_n$ that are divisible by $\rho$.

The following notation will be used in Sections~\ref{proof}
and~\ref{sec:sufficiency}.  For a finite set $\Ac\subseteq S$ and
$\rho\in R_1$, let
\begin{itemize}
    \item $\Ac_\rho=\{\tau\in\Ac\mid\rho\mid\tau\}$;
    \item $\Ac_{\rho^d}=\{\tau\in\Ac\mid m(\tau,\rho)=d\}$
    for each positive multiplicity $d$ occurring in $\Ac_\rho$;
    \item $\Ac_{\rho^d}=\bigsqcup_j\Ac_{\rho^d,j}$ is its
    decomposition into mod $\rho$ equivalence classes.
\end{itemize}
Thus $\Ac_\rho=\bigsqcup_{d,j}\Ac_{\rho^d,j}$.

\begin{definition}[Evenness condition]\label{def:tuple-condition}
Let $\Ac\subseteq S$ be finite and $\rho\in R_1$.  We say that
$\Ac$ satisfies the \emph{evenness condition with respect to $\rho$}
if, for every mod $\rho$ equivalence class $\Ac_{\rho^d,j}$,
\[
 m(\Ac_{\rho^d,j},\gamma)\equiv0\pmod2
 \qquad\text{for every }\gamma\in S_{<d}.
\]
The set $\Ac$ satisfies the \emph{evenness condition} if it does so
with respect to every $\rho\in R_1$.
\end{definition}

\begin{lemma}\label{lem:tuple-test-reduction}
Let $\Ac_{\rho^d,j}$ be a mod $\rho$ equivalence class.  Then
$\Ac_{\rho^d,j}$ satisfies the evenness condition with respect to
$\rho$ if and only if
\[
 m(\Ac_{\rho^d,j},\gamma)\equiv0\pmod2
\]
for those $\gamma\in S_{<d}$ such that $\rho\nmid\gamma$ and
$\gamma\mid\tau$ for at least one $\tau\in\Ac_{\rho^d,j}$.
\end{lemma}

\begin{proof}
The forward implication follows directly from
Definition~\ref{def:tuple-condition}.  Conversely, suppose that the
displayed congruence holds for all the stated tuples.

First let $\gamma\in S_{<d}$ be divisible by $\rho$, and write
$\gamma=\rho^t\gamma_0$, where $t\geq1$ and
$\rho\nmid\gamma_0$.  Since every element of the class contains
$\rho$ with multiplicity $d$, Definition~\ref{multiplicity} gives
\[
 m(\Ac_{\rho^d,j},\gamma)
 =\binom{d}{t}m(\Ac_{\rho^d,j},\gamma_0).
\]
Moreover, $\deg\gamma_0<\deg\gamma<d$.  If $\gamma_0$ divides an
element of the class, its total multiplicity is even by hypothesis;
if it divides none, its total multiplicity is zero.  In either case,
the displayed identity shows that the total multiplicity of $\gamma$
is even.

It remains to consider $\gamma\in S_{<d}$ with
$\rho\nmid\gamma$.  The hypothesis applies if $\gamma$ divides an
element of the class.  Otherwise, for
each $\tau$ there is a $\beta\in R_1$ such that
$m(\gamma,\beta)>m(\tau,\beta)$.  The corresponding binomial
coefficient in Definition~\ref{multiplicity} is zero, and therefore
$m(\tau,\gamma)=0$.  Summing over the class gives
$m(\Ac_{\rho^d,j},\gamma)=0$.
\end{proof}

\begin{example}
    Consider the standard $G_2$-action on $(\mathbb{R}P^2)^3$, where
    $G_2=\mathbb Z_2^2$.  Its fixed point polynomial in the image of $\phi_*$ is
    \[
      (\rho_1\rho_2+\rho_1\rho_{12}+\rho_2\rho_{12})^3,
    \]
    which reduces over $\mathbb Z_2$ to the following nine tangent
    representation classes:
    \[
    \begin{aligned}
    &\rho_1^3\rho_2^3,\ \rho_1^3\rho_{12}^3,\ \rho_2^3\rho_{12}^3,\ 
    \rho_1^3\rho_2^2\rho_{12},\ \rho_1^3\rho_2\rho_{12}^2,\\
    &\rho_1^2\rho_2^3\rho_{12},\ \rho_1^2\rho_2\rho_{12}^3,\ 
    \rho_1\rho_2^2\rho_{12}^3,\ \rho_1\rho_2^3\rho_{12}^2,
    \end{aligned}
    \]
    where $\{\rho_1,\rho_2\}$ is a basis of $\Hom(\Z_2^2,\Z_2)$ and $\rho_{12}=\rho_1+\rho_2$. 
    Denote the resulting set by $\Ac$ and take $\rho=\rho_1$.  Then
    \[ \Ac_{\rho_1} = \Ac \setminus \{ \rho_2^3\rho_{12}^3 \} \]
    consists of the eight weight monomials containing $\rho_1$.  They form
    three mod $\rho_1$ equivalence classes, distinguished by the
    multiplicity of $\rho_1$:
    \[
    \begin{aligned}
    \Ac_{\rho_1^3,1} &= \{ \rho_1^3\rho_2^3,\ \rho_1^3\rho_{12}^3,\ \rho_1^3\rho_2^2\rho_{12},\ \rho_1^3\rho_2\rho_{12}^2 \}, \\
    \Ac_{\rho_1^2,1} &= \{ \rho_1^2\rho_2^3\rho_{12},\ \rho_1^2\rho_2\rho_{12}^3 \}, \\
    \Ac_{\rho_1^1,1} &= \{ \rho_1\rho_2^2\rho_{12}^3,\ \rho_1\rho_2^3\rho_{12}^2 \}.
    \end{aligned}
    \]
    In each class $\Ac_{\rho_1^d,j}$ every weight monomial contains
    $\rho_1$ with the same multiplicity $d$.  
    By Lemma~\ref{lem:tuple-test-reduction}, it suffices to test the
    tuples not divisible by $\rho_1$.  For $\Ac_{\rho_1^3,1}$ these
    have degree at most $2$.
    \begin{itemize}
        \item  $\deg \gamma = 0$, $\gamma=1$: $m(\Ac_{\rho_1^3,1},1)=|\Ac_{\rho_1^3,1}| = 4 $.
        \item  $\deg \gamma = 1$: $\gamma\in \{\rho_2,\rho_{12}\}$ and $m(\Ac_{\rho_1^3,1},\rho_2)=\binom{3}{1}+\binom{0}{1}+\binom{2}{1}+\binom{1}{1} = 6 $, \ $m(\Ac_{\rho_1^3,1},\rho_{12}) = \binom{0}{1}+\binom{3}{1}+\binom{1}{1}+\binom{2}{1} = 6$.
        \item  $\deg \gamma = 2$: $\gamma\in \{\rho_2^2,\rho_{12}^2,\rho_2\rho_{12}\}$ and $m(\Ac_{\rho_1^3,1},\rho_{2}^2) = \binom{3}{2} + \binom{0}{2} + \binom{2}{2} + \binom{1}{2} = 4$, \ $m(\Ac_{\rho_1^3,1},\rho_{12}^2) =\binom{0}{2} + \binom{3}{2} + \binom{1}{2} + \binom{2}{2} = 4 $, \ $m(\Ac_{\rho_1^3,1},\rho_2\rho_{12}) = \binom{3}{1} \binom{0}{1} + \binom{0}{1}\binom{3}{1} + \binom{2}{1}\binom{1}{1} + \binom{1}{1}\binom{2}{1} = 4 $.
    \end{itemize}
    
    For $\Ac_{\rho_1^2,1}$, it suffices to test
    $\gamma\in S_{\rho_1}$ of degree at most $1$.
    \begin{itemize}
        \item  $\deg \gamma = 0$: $m(\Ac_{\rho_1^2,1},1)=|\Ac_{\rho_1^2,1}| = 2 \equiv 0$.
        \item  $\deg \gamma = 1$: $\gamma\in \{\rho_2,\rho_{12}\}$ and $m(\Ac_{\rho_1^2,1},\rho_2)= \binom{3}{1}+\binom{1}{1} = 4 $, \ $m(\Ac_{\rho_1^2,1},\rho_{12}) = \binom{1}{1}+\binom{3}{1} = 4 $.
    \end{itemize}
    
    For $\Ac_{\rho_1^1,1}$,  the evenness condition requires $ m(\Ac_{\rho_1^1,1}, 1) =|\Ac_{\rho_1^1,1}|=2 \equiv 0 \pmod{2}$ for the $0$-tuple.
    
    Thus $\Ac_{\rho_1}$ satisfies the evenness condition with respect
    to $\rho_1$.  The cases $\rho_2$ and $\rho_{12}$ are analogous, so
    $\Ac$ satisfies the evenness condition.
    \end{example}

\begin{remark}[Multiplicity bounds]\label{rem:multiplicity-bounds}
Let $\Ac\subseteq S_n$ be a finite set of faithful weight monomials satisfying
the evenness condition.  If $\tau\in\Ac$, $\rho\mid\tau$, and
$d=m(\tau,\rho)$, then
\[
 d\leq\frac n2
 \qquad\text{and}\qquad
 d\leq n-k+1.
\]
Indeed, if $d>n/2$, let $\Ac_{\rho^d, 1}$ be the mod $\rho$ equivalence class of
$\tau$ and put $g_\tau=\tau/\rho^d$.  Then
$\deg g_\tau=n-d<d$, so the evenness condition applies to
$g_\tau$.  But $g_\tau$ divides no element of $\Ac_{\rho^d, 1}$ except $\tau$
itself, hence $m(\Ac_{\rho^d, 1},g_\tau)=1$, a contradiction.

The second inequality follows from faithfulness. After the $d$ copies
of $\rho$, at least $k-1$ further weights are needed to span
$\Hom(G_k,\mathbb Z_2)$.  Thus $d+k-1\leq n$.

In particular, when $n=k+1$, every repeated weight has multiplicity
at most $2$.  This explains the restriction on repeated weights used
in the computation of $\mathcal Z_4(G_3)$ below.
\end{remark}

\subsection{Regular labeled graphs of $G_k$-manifolds with isolated fixed points}
\label{regular graph}
This subsection gives a geometric interpretation of the empty-tuple
part of the evenness condition of
Definition~\ref{def:tuple-condition}, namely the congruence
\[
 m(\Ac_{\rho^d,j},1)=| \Ac_{\rho^d,j}| \equiv0\pmod2,
\]
i.e.\ the evenness of the cardinality of every mod $\rho$ equivalence
class, and then constructs the \emph{regular labeled graphs} of $M$
used in Section~\ref{sec:z4-application} to visualize the fixed point
data of the explicit manifolds.

Let $M^n$ be a smooth connected closed $G_k$-manifold with isolated fixed
points.  For each $p\in M^{G_k}$, the tangent representation
$\tau_p=T_pM$ is encoded by a faithful weight monomial of degree $n$.

\smallskip
\noindent\emph{The geometric interpretation.}
Fix $\rho\in R_1$ and put $H=\ker\rho$.  Let $C$ be a connected
component of $M^H$ containing a $G_k$-fixed point, and write
$F_C=C\cap M^{G_k}$.  The quotient group $G_k/H\cong\Z_2$ acts on the
closed manifold $C$ with fixed point set
$F_C=C^{G_k/H}$, so by the Conner--Floyd description of
$\mathcal Z_*(\Z_2)$ the set $F_C$ has even cardinality.

For $p\in F_C$, the tangent space $T_pC$ is the $H$-fixed subspace of
$T_pM$.  Since the only weights restricting trivially to $H$ are $0$
and $\rho$, and since $p$ is isolated in $M^{G_k}$, $T_pC$ is the direct sum of $d$ copies of $\rho$ with $d=\dim C$. Hence
$m(\tau_p,\rho)=d$ for every $p\in F_C$.  Moreover, the normal
representations of $C$ in $M$ at two fixed points of the same
component have the same restriction to $H$; indeed, the
$H$-isotypical pieces of the equivariant normal bundle of $C$ have
locally constant ranks, hence constant ranks on the connected manifold
$C$.  Equivalently, their reduced weight monomials are congruent
modulo $\rho$, i.e.\
$\tau_p/\rho^d\equiv\tau_q/\rho^d\pmod\rho$ for all $p,q\in F_C$.
Thus all points of $F_C$ determine tangent representations lying in a
single mod $\rho$ equivalence class.

The sets $F_C$, as $C$ ranges over the components of $M^H$ containing
a $G_k$-fixed point, partition $M^{G_k}$.  Hence every mod $\rho$
equivalence class $\Ac_{\rho^d,j}$ is the disjoint union of the sets
$F_C$ of those components $C$ whose reduced weight monomial $\tau_p/\rho^d$
has the common residue modulo $\rho$, where several components may
contribute to the same class.  Since each $F_C$ has even cardinality,
counting modulo $2$ gives
\[
 m(\Ac_{\rho^d,j},1)=|\Ac_{\rho^d,j}|
 \equiv\sum_C|F_C|\equiv0\pmod2.
\]
 This is exactly
the empty-tuple part of the evenness condition.

\smallskip
\noindent\emph{The construction of the graph.}
For $\rho\in R_1$ and a component $C$ of $M^{\ker\rho}$ containing a
$G_k$-fixed point, the geometric interpretation above gives $|F_C|$
even, so with $d=\dim C$ we may choose a $d$-regular graph, allowing
multiple edges, $\Gamma_{\rho,C}$ with vertex set $F_C$.  This choice
is not unique, but only the resulting local parity information will be
used.  When $d=1$, $C$ is a circle with the standard involution and
$\Gamma_{\rho,C}$ is the unique edge joining its two fixed points.

Taking the union of all graphs $\Gamma_{\rho,C}$, as $\rho$ and
$C$ vary, gives an $n$-regular graph $\Gamma_M$ with vertex set
$M^{G_k}$.  We label each edge by the weight corresponding to its
tangent direction:
\[
 \alpha\colon E_{\Gamma_M}\longrightarrow R_1 .
\]
The pair $(\Gamma_M,\alpha)$ is called a \emph{regular labeled graph}
of $M$ \cite{L3}.  It satisfies the following two properties.
\begin{enumerate}[label=(\arabic*)]
\item For each vertex $p$, the product of the labels on the edges
incident to $p$ is exactly $\tau_p$.  In particular, these labels span
$\Hom(G_k,\Z_2)$.
\item For each subgraph $\Gamma_{\rho,C}$ and any two vertices
$p,q\in V_{\Gamma_{\rho,C}}$,
\[
 \tau_p/\rho^d\equiv \tau_q/\rho^d\pmod\rho,
\]
and $|V_{\Gamma_{\rho,C}}|$ is even.
\end{enumerate}

Property (1) holds by the definition of the labels and the
effectiveness of the action.  Property (2) restates the mod $\rho$
congruence $\tau_p/\rho^d\equiv\tau_q/\rho^d\pmod\rho$ and the
evenness of $|F_C|$ established in the geometric interpretation
above.

\section{Necessity of the Theorem~\ref{A}}\label{proof}

\subsection{The local parity criterion}

The necessity direction has two ingredients.  The first is a
tuple-shift lemma: once the lower-degree multiplicities vanish, the
multiplicity of a tuple is determined, modulo $2$, by its reduction
modulo $\rho$.  The second shows that the localized integrality coming
from a normal bundle forces exactly these lower-degree vanishings.

\begin{lemma}\label{lem:tuple-congruence}
Let $\Ac_{\rho^d,1}\subseteq S_n$ be a mod $\rho$
equivalence class, where $m(\tau,\rho)=d>0$ for every $\tau$.
Extend $\rho$ to a basis $\{\rho,\rho_2,\dots,\rho_k\}$ of
$\Hom(G_k,\mathbb Z_2)$ and put
$R_0=\mathbb Z_2[\rho_2,\dots,\rho_k]$.  Let $g_0\in R_0$ be the
common reduction modulo $\rho$ of $\tau/\rho^d$.

Suppose that
\[
 m(\Ac_{\rho^d,1},\delta)\equiv0\pmod2
 \quad\text{whenever }\rho\nmid\delta
 \text{ and }\deg\delta\leq d-2.
\]
Let
$\gamma=\prod_{i=1}^s\beta_i^{u_i}$ divide $g_0$, where the
$\beta_i\in R_0$ are distinct and $\sum_i u_i=d-1$.  Then, for
$0\leq b_i\leq u_i$ and
\[
 \gamma'=\prod_{i=1}^s
 \beta_i^{u_i-b_i}(\beta_i+\rho)^{b_i},
\]
one has
\begin{equation}\label{eq:lemma-target}
m(\Ac_{\rho^d,1},\gamma')
\equiv
\big(\prod_{i=1}^s\binom{u_i}{b_i}\big)\,
m(\Ac_{\rho^d,1},\gamma)\pmod2.
\end{equation}
\end{lemma}

\begin{proof}
We prove the assertion by induction on
\[
 \ell=|\{i\mid b_i>0\}|.
\]
The case $\ell=0$ is immediate. Suppose that $\ell>0$ and that the
assertion holds whenever fewer than $\ell$ of the integers $b_i$ are
positive. After reindexing, assume that $b_s>0$.
 Set
\[
 \gamma_0=\prod_{i<s}\beta_i^{u_i},
 \qquad
 \gamma'_0=\prod_{i<s}
 \beta_i^{u_i-b_i}(\beta_i+\rho)^{b_i}.
\]
Let
\[
 N_s=m(g_0,\beta_s),\qquad
 x_{s,\tau}=m(\tau,\beta_s).
\]
Because $g_0$ records the weights modulo $\rho$, the multiplicity of
$\beta_s+\rho$ in $\tau$ is $N_s-x_{s,\tau}$.  Hence
\[
 m(\tau,\gamma')=m(\tau,\gamma'_0)m(\tau, \beta_s^{u_s-b_s}( \beta_s+\rho)^{b_s})=m(\tau,\gamma'_0)
 \binom{N_s-x_{s,\tau}}{b_s}
 \binom{x_{s,\tau}}{u_s-b_s}.
\]
The left-hand side below is an integer-valued polynomial of degree
$u_s$ in $x$.  It therefore has a unique expansion in the binomial
basis, 
\[
 \binom{N_s-x}{b_s}\binom{x}{u_s-b_s}
 =\sum_{v=0}^{u_s}c_v\binom{x}{v},
 \qquad c_v\in\mathbb Z.
\]
The coefficient of $x^{u_s}$ on the left is
$(-1)^{b_s}/(b_s!(u_s-b_s)!)$, whereas the leading coefficient of
$\binom{x}{u_s}$ is $1/u_s!$.  Hence
\[
 c_{u_s}=(-1)^{b_s}\binom{u_s}{b_s}
 \equiv\binom{u_s}{b_s}\pmod2.
\]
Substituting $x=x_{s,\tau}$ into the binomial-basis expansion,
multiplying by $m(\tau,\gamma'_0)$, and summing over the class gives
\begin{align*}
 m(\Ac_{\rho^d,1},\gamma')
 &=
 \sum_{v=0}^{u_s}c_v
 \sum_{\tau\in\Ac_{\rho^d,1}}
 m(\tau,\gamma'_0)\binom{x_{s,\tau}}v\\
 &=
 \sum_{v=0}^{u_s}c_v\,
 m(\Ac_{\rho^d,1},\gamma'_0\beta_s^v).
\end{align*}
For $v<u_s$,
 $\gamma'_0\beta_s^v$ is not divisible by $\rho$ and has degree at
most $d-2$.  The parity hypothesis therefore eliminates all such
terms.  It follows that
\[
 m(\Ac_{\rho^d,1},\gamma')\equiv c_{u_s}m(\Ac_{\rho^d,1},\gamma'_0\beta_s^{u_s})
 \equiv \binom{u_s}{b_s}
 m(\Ac_{\rho^d,1},\gamma'_0\beta_s^{u_s})
 \pmod2.
\]
The tuple $\gamma'_0\beta_s^{u_s}$ has $\ell-1$ shifted components.
The induction hypothesis gives
\[
 m(\Ac_{\rho^d,1},\gamma'_0\beta_s^{u_s})
 \equiv
\big( \prod_{i=1}^{s-1}\binom{u_i}{b_i}\big)\,
 m(\Ac_{\rho^d,1},\gamma)\pmod2.
\]
Combining this congruence with the preceding one proves
\eqref{eq:lemma-target}.
\end{proof}

\begin{lemma}[Normal-bundle integrality implies evenness]
\label{lem:normal-integrality-implies-evenness}
Let $\Ac_{\rho^d,1}\subseteq S_n$ be a mod $\rho$ equivalence
class with $m(\tau,\rho)=d>0$ for every
$\tau\in\Ac_{\rho^d,1}$, and put
$g_\tau=\tau/\rho^d$.  Suppose that
\[
 \frac{1}{\rho^d}
 \sum_{\tau\in\Ac_{\rho^d,1}}h(g_\tau)
 \in S_\rho^{-1}R
\]
for every symmetric polynomial $h$.  Then the class $\Ac_{\rho^d,1}$ satisfies the
evenness condition with respect to $\rho$.
\end{lemma}

\begin{proof}
Extend $\rho$ to a basis $\{\rho,\rho_2,\dots,\rho_k\}$ of
$\Hom(G_k,\mathbb Z_2)$, and set
$R_0=\mathbb Z_2[\rho_2,\dots,\rho_k]$.  The mod $\rho$ equivalence
class has a unique weight monomial $g_0\in R_0$ such that
\[
 g_\tau\equiv g_0\pmod\rho
 \qquad \forall \tau\in\Ac_{\rho^d,1}.
\]
Indeed, reducing a weight different from $\rho$ modulo $\rho$ produces a
nonzero linear form in $R_0$, and the defining congruence of the class
shows that the resulting product is independent of $\tau$.

For a fixed $h$, write
$F_h=\sum_{\tau\in\Ac_{\rho^d,1}}h(g_\tau)\in R$.  The condition
$F_h/\rho^d\in S_\rho^{-1}R$ means that there is an
$s\in S_\rho$ for which $\rho^d$ divides $sF_h$.  Since
$\rho\nmid s$ and $R$ is a unique factorization domain, this is
equivalent to $\rho^d\mid F_h$.  Thus the hypothesis is equivalent
to
\begin{equation}\label{eq:integral-condition-2-refined}
\sum_{\tau\in\Ac_{\rho^d,1}}h(g_\tau)
\equiv0\pmod{\rho^d}
\end{equation}
for every symmetric polynomial $h$.  We argue by induction on $d$.

\smallskip
\noindent\emph{Base case.}
If $d=1$, take $h=1$.  The hypothesis gives
$|\Ac_{\rho^1,1}|/\rho\in S_\rho^{-1}R$, so $\rho\mid|\Ac_{\rho^1,1}|$
in $R$ and hence $|\Ac_{\rho,1}|$ is even.  This is the evenness
condition when $d=1$.

\smallskip
\noindent\emph{Inductive step.}
Assume $d>1$ and that the result holds for $d-1$.  The set
\[
 \mathcal B=\{\tau/\rho\mid\tau\in\Ac_{\rho^d,1}\}
\]
is a mod $\rho$ equivalence class with $\rho$-multiplicity $d-1$.
Moreover,
\[
 \frac{1}{\rho^{d-1}}\sum_{\tau'\in\mathcal B}
 h\left(\frac{\tau'}{\rho^{d-1}}\right)
 =
 \rho\left(
 \frac{1}{\rho^d}\sum_{\tau\in\Ac_{\rho^d,1}}h(g_\tau)
 \right)
 \in S_\rho^{-1}R.
\]
The induction hypothesis therefore applies to $\mathcal B$.  For a tuple $\delta$, if
$\rho\nmid\delta$ and $\deg\delta\leq d-2$, then
$m(\tau,\delta)=m(\tau/\rho,\delta)$, and hence
\[
 m(\Ac_{\rho^d,1},\delta)=m(\mathcal B,\delta)
 \equiv0\pmod2.
\]
By Lemma~\ref{lem:tuple-test-reduction}, it remains only to prove this
congruence for the $(d-1)$-tuples not divisible by $\rho$ that divide
some member of the class.

Put $q=n-d=\deg g_\tau$.  If $q<d-1$, no such tuple can divide
$g_\tau$, and the proof is complete.  We may therefore assume
$q\geq d-1$.

For an odd positive integer $r$, take
\[
 h_r(x_1,\dots,x_q)
 =
 \sum_{1\leq i_1<\cdots<i_{d-1}\leq q}
 (x_{i_1}\cdots x_{i_{d-1}})^r.
\]
Then
\[
 h_r(g_\tau)=
 \sum_{\substack{\gamma'\mid g_\tau\\
                  \deg\gamma'=d-1}}
 m(\tau,\gamma')(\gamma')^r.
\]
Indeed, the definition of $h_r$ sums over all choices of
$d-1$ weights of $g_\tau$, and exactly
$m(\tau,\gamma')$ of these choices have product $\gamma'$.
Every tuple $\gamma'$ occurring in this sum has a unique unshifted
reduction modulo $\rho$.  If that reduction is
\[
 \gamma=\prod_{i=1}^s\beta_i^{u_i}\mid g_0,
\]
where $u_i>0$ and the $\beta_i\in R_0$ are distinct, then $\gamma'$ has the form
\[
 \gamma_{\mathbf b}
 :=
 \prod_{i=1}^s
 \beta_i^{u_i-b_i}(\beta_i+\rho)^{b_i},
 \qquad
 \mathbf b=(b_1,\dots,b_s),\quad 0\leq b_i\leq u_i.
\]
Conversely, a tuple of this form that does not divide any $g_\tau$
has total multiplicity zero.  We may therefore group the preceding
sum by $\gamma$ and write
\begin{equation}\label{eq:integral-condition-3-refined}
\sum_{\tau\in\Ac_{\rho^d,1}}h_r(g_\tau)
=
\sum_{\substack{\gamma\mid g_0\\ \deg\gamma=d-1}}
\ \sum_{\substack{0\leq b_i\leq u_i\\1\leq i\leq s}}
m(\Ac_{\rho^d,1},\gamma_{\mathbf b})
(\gamma_{\mathbf b})^r
\equiv0\pmod{\rho^d},
\end{equation}
where the outer sum runs over the distinct $(d-1)$-tuples dividing
$g_0$.

Fix $\gamma=\prod_{i=1}^s\beta_i^{u_i}$ in the outer sum, and denote
its inner sum by $S_\gamma$.  By the definition of
$\gamma_{\mathbf b}$,
\begin{align*}
S_\gamma
&=
\sum_{b_1=0}^{u_1}\cdots\sum_{b_s=0}^{u_s}
m(\Ac_{\rho^d,1},\gamma_{\mathbf b})
\prod_{i=1}^s
\beta_i^{(u_i-b_i)r}(\beta_i+\rho)^{b_i r}.
\end{align*}
Lemma~\ref{lem:tuple-congruence} gives
\[
 m(\Ac_{\rho^d,1},\gamma_{\mathbf b})
 \equiv
 \left(\prod_{i=1}^s\binom{u_i}{b_i}\right)
 m(\Ac_{\rho^d,1},\gamma)
 \pmod2.
\]
Since the coefficient ring is $\mathbb Z_2$, this congruence is an
equality in $R$.  Substituting it into the preceding sum and then
separating the sums over $b_1,\dots,b_s$, we obtain
\begin{align*}
S_\gamma
&=
m(\Ac_{\rho^d,1},\gamma)
\sum_{b_1=0}^{u_1}\cdots\sum_{b_s=0}^{u_s}
\prod_{i=1}^s
\binom{u_i}{b_i}
(\beta_i^r)^{u_i-b_i}
\bigl((\beta_i+\rho)^r\bigr)^{b_i}\\
&=
m(\Ac_{\rho^d,1},\gamma)
\prod_{i=1}^s
\left[
\sum_{b_i=0}^{u_i}
\binom{u_i}{b_i}
(\beta_i^r)^{u_i-b_i}
\bigl((\beta_i+\rho)^r\bigr)^{b_i}
\right]\\
&=
m(\Ac_{\rho^d,1},\gamma)
\prod_{i=1}^s
\left(\beta_i^r+(\beta_i+\rho)^r\right)^{u_i},
\end{align*}
where the last equality is the binomial theorem over $\mathbb Z_2$.

We next determine the lowest power of $\rho$ in each factor.  Since
$r$ is odd,
\[
 \beta_i^r+(\beta_i+\rho)^r
 \equiv\rho\beta_i^{r-1}\pmod{\rho^2},
\]
and therefore
\[
 \left(\beta_i^r+(\beta_i+\rho)^r\right)^{u_i}
 \equiv
 \rho^{u_i}\beta_i^{(r-1)u_i}
 \pmod{\rho^{u_i+1}}.
\]
Multiplying these congruences and using
$\sum_i u_i=d-1$ gives
\[
 \prod_{i=1}^s
 \left(\beta_i^r+(\beta_i+\rho)^r\right)^{u_i}
 \equiv
 \rho^{d-1}\prod_{i=1}^s\beta_i^{(r-1)u_i}
 =
 \rho^{d-1}\gamma^{r-1}
 \pmod{\rho^d}.
\]
Substitution in \eqref{eq:integral-condition-3-refined} gives
\[
 \rho^{d-1}
 \sum_{\substack{\gamma\mid g_0\\\deg\gamma=d-1}}
 m(\Ac_{\rho^d,1},\gamma)\gamma^{r-1}
 \equiv0\pmod{\rho^d}.
\]
The sum  belongs to $R_0$.  Since
$R=R_0[\rho]$, the preceding congruence implies that this sum is
divisible by $\rho$.  Its reduction in $R/(\rho)=R_0$ is the sum
itself, so it must be zero.  We have proved
\begin{equation}\label{eq:vandermonde-system}
 \sum_{\substack{\gamma\mid g_0\\\deg\gamma=d-1}}
 m(\Ac_{\rho^d,1},\gamma)\gamma^{r-1}=0
 \quad\text{in }R_0.
\end{equation}

Let $\gamma_1,\dots,\gamma_N$ be the distinct tuples in this sum and
put
\[
 a_j=m(\Ac_{\rho^d,1},\gamma_j)\in\mathbb Z_2.
\]
If $N=0$, there is nothing left to prove.  Assume $N>0$.
For $r=2\ell+1$, equation~\eqref{eq:vandermonde-system} becomes
\[
 \sum_{j=1}^N a_j(\gamma_j^2)^\ell=0.
\]
Taking $\ell=0,\dots,N-1$ gives the linear system
\[
 \begin{pmatrix}
 1&1&\cdots&1\\
 \gamma_1^2&\gamma_2^2&\cdots&\gamma_N^2\\
 \vdots&\vdots&&\vdots\\
 (\gamma_1^2)^{N-1}&(\gamma_2^2)^{N-1}
     &\cdots&(\gamma_N^2)^{N-1}
 \end{pmatrix}
 \begin{pmatrix}a_1\\a_2\\ \vdots\\a_N\end{pmatrix}
 =0.
\]
The determinant of this Vandermonde matrix is
\[
 \prod_{a<b}(\gamma_a^2+\gamma_b^2)
 =\prod_{a<b}(\gamma_a+\gamma_b)^2\neq0.
\]
Indeed, the $\gamma_j$ are distinct elements of the integral domain
$R_0$.
The matrix is therefore invertible over the fraction field of $R_0$,
so $a_1=\cdots=a_N=0$.

Finally, let $\gamma'$ be any $(d-1)$-tuple not divisible by $\rho$
that divides some member of the class.  Its unshifted reduction is one
of the $\gamma_j$, and $\gamma'=(\gamma_j)_{\mathbf b}$ for a suitable
$\mathbf b$.  Lemma~\ref{lem:tuple-congruence} now gives
\[
 m(\Ac_{\rho^d,1},\gamma')
 \equiv
 \left(\prod_i\binom{u_i}{b_i}\right)a_j=0\pmod2.
\]
Together with the low-degree case and
Lemma~\ref{lem:tuple-test-reduction}, this proves the evenness
condition with respect to $\rho$.
\end{proof}

\subsection{Proof of necessity}

\begin{proof}[Proof of necessity in the Theorem~\ref{A}]
Let $M^n$ be a $G_k$-manifold with isolated fixed point set
and its fixed point data
\[
 \Ac=\{\tau_p\mid p\in M^{G_k}\}.
\]
It is enough to verify the evenness condition for each
$\rho\in R_1$ which occurs in the fixed point data.

\smallskip
\noindent\emph{Step 1: the class decomposition.}
Fix such a $\rho$ and put $H=\ker\rho$.  By the geometric
interpretation in Subsection~\ref{regular graph}, the sets
$F_C=C\cap M^{G_k}$ (as $C$ ranges over the components of $M^H$
containing a $G_k$-fixed point) have even cardinality and
$\rho$-multiplicity $d=\dim C$, each $F_C$ lies in a single mod
$\rho$ equivalence class, and every class $\Ac_{\rho^d,j}$ is the
disjoint union of the sets $F_C$ of those components whose reduced
weight monomial has the prescribed residue modulo $\rho$.  

\smallskip
\noindent\emph{Step 2: localization on each class.}
Let $\nu_C=\nu(C,M)$ be the $G_k$-equivariant normal bundle of $C$.
By the equivariant tubular neighborhood theorem, a $G_k$-invariant
neighborhood of $C$ in $M$ is equivariantly diffeomorphic to a
neighborhood of the zero section in $\nu_C$.  Since
$T_pC$ is the direct sum of $d$ copies of $\rho$, the weight monomial
$g_p=\tau_p/\rho^d$ is the product of the normal weights at $p$.  Apply the
Kosniowski--Stong localization formula \cite{KS} to the $G_k$-bundle
$\nu_C$ over $C$, for every symmetric polynomial $h$ in the
Stiefel--Whitney roots of $\nu_C$,
\begin{equation}\label{eq:KS-normal-bundle}
 \frac{1}{\rho^d}
 \sum_{p\in F_C}h(g_p)
 \in R\subseteq S_\rho^{-1}R,
\end{equation}
because the equivariant Euler class of $T_pC$ is $\rho^d$.

Summing \eqref{eq:KS-normal-bundle} over the components $C$ contained
in a class $\Ac_{\rho^d,j}$ (see Step 1) yields
\[
 \frac{1}{\rho^d}
 \sum_{\tau\in\Ac_{\rho^d,j}}h(g_\tau)
 \in S_\rho^{-1}R
\]
for every symmetric polynomial $h$.  Lemma~\ref{lem:normal-integrality-implies-evenness}
therefore shows that
$\Ac_{\rho^d,j}$ satisfies the evenness condition with respect to
$\rho$.

This applies to every mod $\rho$ equivalence class.  If a weight
$\rho\in R_1$ does not occur in the fixed point data, the condition
with respect to $\rho$ is vacuous.  Hence $\Ac$ satisfies the
evenness condition.
\end{proof}

\section{Sufficiency of the Theorem~\ref{A}}\label{sec:sufficiency}

We first record the algebraic facts used to pass from the evenness
condition to localized integrality.  They are elementary, but we keep
them separate so that the proof of sufficiency has a clear structure:
expand the numerator, identify each low-order coefficient with a
tuple multiplicity, and then apply the evenness condition.

\subsection{Two algebraic tools for local integrality}\label{subsec:local-technical}

\subsubsection{Symmetric polynomials under linear shifts}

Let $g_0=\beta_1\cdots\beta_q$ be a weight monomial not divisible by $\rho$.
For an index set $L\subseteq\{1,\dots,q\}$, write
\[
 \beta_L=\prod_{i\in L}\beta_i,\qquad
 \bar\beta_L=\prod_{i\in L}(\beta_i+\rho).
\]

\begin{lemma}\label{lem:expansion-total}
Let $I\sqcup J=\{1,\dots,q\}$ and
$g'=\bar\beta_I\beta_J$.  For $1\leq r\leq q$,
\[
 e_r(g')
 =
 \sum_{u=0}^r\rho^u
 \sum_{\substack{\beta_L\mid g_0\\\deg\beta_L=u}}
 m(g',\bar\beta_L)\,
 e_{r-u}(g_0/\beta_L),
\]
where the inner sum runs over the distinct degree-$u$ subproducts of
$g_0$.
\end{lemma}

\begin{proof}
Expanding each shifted factor $\bar\beta_i=\beta_i+\rho$ gives
\[
 e_r(g')=
 \sum_{\substack{K\subseteq I\sqcup J\\|K|=r}}
 \ \sum_{U\subseteq K\cap I}
 \rho^{|U|}\prod_{s\in K\setminus U}\beta_s.
\]
Here $K$ records the $r$ factors chosen in the elementary symmetric
polynomial, while $U$ records those chosen factors from which the
$\rho$-term is taken.  Fix $U\subseteq I$.  Any set $K$ contributing
to the inner sum can then be written uniquely as
$K=U\sqcup V$, where
\[
 V\subseteq\{1,\dots,q\}\setminus U,
 \qquad |V|=r-|U|.
\]
Therefore, after interchanging the two sums,
\begin{align*}
 e_r(g')
 &=
 \sum_{\substack{U\subseteq I\\|U|\leq r}}
 \rho^{|U|}
 \sum_{\substack{V\subseteq\{1,\dots,q\}\setminus U\\
                  |V|=r-|U|}}
 \prod_{s\in V}\beta_s\\
 &=
 \sum_{\substack{U\subseteq I\\|U|\leq r}}
 \rho^{|U|}e_{r-|U|}(g_0/\beta_U).
\end{align*}
For a fixed degree-$u$ subproduct $\beta_L$ of $g_0$, the number of
subsets $U\subseteq I$ satisfying $\beta_U=\beta_L$ is precisely
$m(g',\bar\beta_L)$: such a subset chooses an occurrence of each
shifted factor in the tuple $\bar\beta_L$ from among the shifted
factors of $g'$.  Grouping the terms with $|U|=u$ and
$\beta_U=\beta_L$ proves the formula.
\end{proof}

\subsubsection{Applications of Lucas's theorem}

By Lucas's theorem~\cite{f}, $\binom{a}{b}$ is odd if and only if every
$1$-bit in the binary expansion of $b$ is also a $1$-bit in that of
$a$.

\begin{lemma}\label{lucas_cor}
Let $\mathcal M=\{\beta_1,\dots,\beta_t\}$ be a multiset of linear
forms, let $\mathcal M_0$ be its underlying set, and let
$v_1,\dots,v_t$ be positive integers.  For each
$\beta\in\mathcal M_0$, let $u_\beta$ be the integer whose $1$-bits
are the union of the $1$-bits of those $v_i$ for which
$\beta_i=\beta$, and set
\[
 \gamma=\prod_{\beta\in\mathcal M_0}\beta^{u_\beta}.
\]
Then
\[
 \gamma\mid\prod_{i=1}^t\beta_i^{v_i},
 \qquad
 \deg\gamma\leq\sum_{i=1}^t v_i,
\]
and, for every $\tau\in S$,
\[
 \prod_{i=1}^t m(\tau,\beta_i^{v_i})
 \equiv m(\tau,\gamma)\pmod2.
\]
\end{lemma}

\begin{proof}
Fix $\tau\in S$.  For $\beta\in\mathcal M_0$, put
$N_\beta=m(\tau,\beta)$.  Lucas's
theorem gives
\[
 \prod_{i:\,\beta_i=\beta}\binom{N_\beta}{v_i}
 \equiv \binom{N_\beta}{u_\beta}\pmod2,
\]
because either side is odd if and only if every $1$-bit occurring in
one of the relevant $v_i$ also occurs in $N_\beta$.  Multiplying over
$\beta\in\mathcal M_0$ yields
\begin{align*}
 \prod_{i=1}^t m(\tau,\beta_i^{v_i})
 &=
 \prod_{\beta\in\mathcal M_0}
 \prod_{i:\,\beta_i=\beta}
 \binom{N_\beta}{v_i}\\
 &\equiv
 \prod_{\beta\in\mathcal M_0}\binom{N_\beta}{u_\beta}
 =m(\tau,\gamma)\pmod2.
\end{align*}
Finally,
$u_\beta\leq\sum_{i:\,\beta_i=\beta}v_i$ for every $\beta$.
This proves the divisibility assertion and the degree bound.
\end{proof}

\subsection{The local integrality criterion}

\begin{lemma}[Evenness implies local integrality]
\label{lem:evenness-implies-local-integrality}
Let $\Ac_{\rho^d,1}\subseteq S_n$ be a mod $\rho$ equivalence
class with $m(\tau,\rho)=d>0$ for every $\tau$.  If this class
satisfies the evenness condition with respect to $\rho$, then
\[
 \sum_{\tau\in\Ac_{\rho^d,1}}\frac{f(\tau)}{\tau}
 \in S_\rho^{-1}R
\]
for every symmetric polynomial $f$ in $n$ variables.
\end{lemma}

\begin{proof}
Extend $\rho$ to a basis $\{\rho,\rho_2,\dots,\rho_k\}$ of
$\Hom(G_k,\mathbb Z_2)$ and put
$R_0=\mathbb Z_2[\rho_2,\dots,\rho_k]$.  Set $q=n-d$ and, for each
$\tau$ in the class, put $g_\tau=\tau/\rho^d$.  There is a unique
weight monomial
\[
 g_0=\beta_1\cdots\beta_q\in R_0
\]
such that $g_\tau\equiv g_0\pmod\rho$ for every $\tau$ in the class.
For each $\tau$, choose a partition
$I_\tau\sqcup J_\tau=\{1,\dots,q\}$ such that
\[
 g_\tau=\beta_{I_\tau}\bar\beta_{J_\tau},
 \qquad
 \bar{g_\tau}=\bar\beta_{I_\tau}\beta_{J_\tau}.
\]
Set
\[
 \bar{g_0}=\prod_{i=1}^q(\beta_i+\rho).
\]
Then $g_\tau\bar{g_\tau}=g_0\bar{g_0}$, and hence
\begin{equation}\label{eq:common-denominator}
 \sum_\tau\frac{h(g_\tau)}{g_\tau}
 =
 \frac{1}{g_0\bar{g_0}}
 \sum_\tau \bar{g_\tau} h(g_\tau).
\end{equation}
Both $g_0$ and $\bar{g_0}$ belong to $S_\rho$.

\smallskip
\noindent\emph{Step 1:}
We first show that
\[
 \frac{1}{\rho^d}
 \sum_{\tau\in\Ac_{\rho^d,1}}\frac{h(g_\tau)}{g_\tau}
 \in S_\rho^{-1}R
\]
for every symmetric polynomial $h$ in $q$ variables.  By
\eqref{eq:common-denominator}, it is enough to show that the numerator
$\sum_\tau \bar{g_\tau} h(g_\tau)$ is divisible by $\rho^d$.  We check
the coefficients of $\rho^s$ for $s<d$.  By linearity, it is enough
to take
$h=e_\lambda=e_{\lambda_1}\cdots e_{\lambda_t}$, We omit the
vanishing cases $\lambda_i>q$.  Applying
Lemma~\ref{lem:expansion-total} gives
\begin{align*}
 e_{\lambda_i}(g_\tau)
 &=
 \sum_{u_i=0}^{\lambda_i}\rho^{u_i}
 \sum_{\substack{\beta_{L_i}\mid g_0\\|L_i|=u_i}}
 m(\tau,\bar\beta_{L_i})\,
 e_{\lambda_i-u_i}(g_0/\beta_{L_i}),\\
 \bar{g_\tau}
 &=
 \sum_{v=0}^{q}\rho^v
 \sum_{\substack{\beta_G\mid g_0\\|G|=v}}
 m(\tau,\beta_G)\,\frac{g_0}{\beta_G}.
\end{align*}
In the second identity, a shifted factor of $\bar{g_\tau}$ corresponds
to an unshifted factor of $g_\tau$.

Multiplying the preceding expansions gives
\begin{align}
 \bar{g_\tau} h(g_\tau)
 ={}&
 \sum_{v=0}^q
 \sum_{u_1=0}^{\lambda_1}\cdots
 \sum_{u_t=0}^{\lambda_t}
 \rho^{v+u_1+\cdots+u_t}
 \notag\\[-2mm]
 &\cdot
 \sum_{\substack{\beta_G\mid g_0\\\deg\beta_G=v}}
 \ \sum_{\substack{\beta_{L_1}\mid g_0\\\deg\beta_{L_1}=u_1}}
 \cdots
 \sum_{\substack{\beta_{L_t}\mid g_0\\\deg\beta_{L_t}=u_t}}
 \frac{g_0}{\beta_G}
 \left(\prod_{i=1}^t
 e_{\lambda_i-u_i}(g_0/\beta_{L_i})\right)
 m(\tau,\beta_G)
 \prod_{i=1}^t m(\tau,\bar\beta_{L_i}).
 \label{eq:expanded-local-numerator}
\end{align}
All inner sums run over distinct subproducts of $g_0$.

Fix one summand in \eqref{eq:expanded-local-numerator}.  Its factor
in $R_0$ is independent of $\tau$; the only part depending on $\tau$
is
\[
 m(\tau,\beta_G)
 \prod_{i=1}^t m(\tau,\bar\beta_{L_i}).
\]
Applying Lemma~\ref{lucas_cor} to the weights occurring in
$\beta_G,\bar\beta_{L_1},\dots,\bar\beta_{L_t}$ gives a tuple
$\gamma$, independent of $\tau$, such that
\begin{equation}\label{eq:lucas-local-coefficient}
 m(\tau,\beta_G)
 \prod_{i=1}^t m(\tau,\bar\beta_{L_i})
 \equiv m(\tau,\gamma)\pmod2,
\end{equation}
and
\[
 \rho\nmid\gamma,
 \qquad
 \deg\gamma
 \leq v+u_1+\cdots+u_t.
\]
The first assertion holds because none of the weights occurring in
$\beta_G$ or $\bar\beta_{L_i}$ is $\rho$.  If all these tuples are
empty, we use the convention $\gamma=1$; the same congruence and
degree bound then hold.

Now the sum \eqref{eq:expanded-local-numerator} over the class and fix
$s<d$.  Every summand in the coefficient of $\rho^s$ has
$v+u_1+\cdots+u_t=s$.  After summing over $\tau$, the
$\tau$-dependent part of that summand becomes
$m(\Ac_{\rho^d,1},\gamma)$ by
\eqref{eq:lucas-local-coefficient}.  Since
$\rho\nmid\gamma$ and $\deg\gamma\leq s<d$, the evenness condition
gives
\[
 m(\Ac_{\rho^d,1},\gamma)=0\quad\text{in }\mathbb Z_2.
\]
Thus every summand in the coefficient of $\rho^s$ is zero.
Hence the coefficients of $\rho^s$ vanish for all $s<d$, and
\[
 \sum_{\tau\in\Ac_{\rho^d,1}}\bar{g_\tau} h(g_\tau)
 \in\rho^d R.
\]
Equation~\eqref{eq:common-denominator} therefore implies
\[
 \frac{1}{\rho^d}
 \sum_{\tau\in\Ac_{\rho^d,1}}\frac{h(g_\tau)}{g_\tau}
 \in S_\rho^{-1}R
\]
for every symmetric polynomial $h$ in $q$ variables.  This completes
Step 1.

\smallskip
\noindent\emph{Step 2: integrality of $\sum_\tau f(\tau)/\tau$.}
Let $f$ be a symmetric polynomial in $n$ variables.  Since
$\tau=\rho^d g_\tau$, substituting the $d$ copies of the weight $\rho$
and the $q$ weights of $g_\tau$ into $f$ gives
\[
 f(\tau)=\sum_a \rho^a h_a(g_\tau),
\]
where each $h_a$ is symmetric in the weights of $g_\tau$.  Applying
Step 1 to each $h_a$ yields
\[
 \sum_{\tau\in\Ac_{\rho^d,1}}\frac{f(\tau)}{\tau}
 =
 \sum_a \rho^a
 \left(
 \frac{1}{\rho^d}
 \sum_{\tau\in\Ac_{\rho^d,1}}
 \frac{h_a(g_\tau)}{g_\tau}
 \right)
 \in S_\rho^{-1}R.
\]
This is the required local integrality.
\end{proof}

\subsection{Proof of sufficiency}

\begin{proof}[Proof of sufficiency in the Theorem~\ref{A}]
Suppose that $\Ac$ satisfies the evenness condition.  Fix
$\rho\in R_1$ and partition
\[
 \Ac_\rho=\bigsqcup_{d,j}\Ac_{\rho^d,j}
\]
into its mod $\rho$ equivalence classes.  Lemma~\ref{lem:evenness-implies-local-integrality}, applied to each class, gives
\[
 \sum_{\tau\in\Ac_{\rho^d,j}}\frac{f(\tau)}{\tau}
 \in S_\rho^{-1}R
\]
for every symmetric polynomial $f$.  Summing over $d$ and $j$ proves
that
\[
 \sum_{\tau\in\Ac_\rho}\frac{f(\tau)}{\tau}
 \in S_\rho^{-1}R .
\]
If $\rho\nmid\tau$, then $\tau\in S_\rho$ and
$f(\tau)/\tau\in S_\rho^{-1}R$.  Hence the full localization sum
\[
 \sum_{\tau\in\Ac}\frac{f(\tau)}{\tau}
\]
belongs to $S_\rho^{-1}R$.  Since this holds for every
$\rho\in R_1$, Lemma~\ref{cap} implies that the sum actually lies in
$R$.  The tom Dieck-Kosniowski--Stong criterion, Theorem~\ref{Kos-Stong}, then
shows that $\Ac$ is the fixed point data of an $n$-dimensional
$G_k$-manifold.
This proves the sufficiency direction and completes the proof of the
Theorem~\ref{A}.
\end{proof}

\section{Application on $\mathcal Z_4(G_3)$}
\label{sec:z4-application}

Let $G_3=(\mathbb Z_2)^3$.  Choose a basis
$\rho_1,\rho_2,\rho_3$ of $\Hom(G_3,\mathbb Z_2)$, and write
\[
 \rho_{ij}=\rho_i+\rho_j,\qquad
 \rho_{123}=\rho_1+\rho_2+\rho_3.
\]
Thus the seven nonzero weights are
\[
\rho_1,\rho_2,\rho_3,\rho_{12},\rho_{13},\rho_{23},\rho_{123}.
\]

This section applies the Theorem~\ref{A} to the first non-square-free
case, namely degree-four fixed point data for $G_3$.  We first reduce
the calculation to the three possible types of faithful monomials and
determine the dimension of the resulting space.  We then realize the
corresponding generators by explicit equivariant manifolds.

\subsection{The dimension count}

Let $E=\Hom(G_3,\mathbb Z_2)$ and
$\Gamma=\operatorname{Aut}(G_3)\cong\operatorname{GL}(E)$.  For $\sigma\in\Gamma$, precomposing a
$G_3$-action $a$ with $\sigma^{-1}$ gives
$(\sigma\cdot a)(g,x)=a(\sigma^{-1}g,x)$ and replaces each tangent
weight $\rho$ by $\rho\circ\sigma^{-1}$.  This defines the
$\Gamma$-action on fixed-point data.
There are two cases.  If a weight is repeated, faithfulness forces the
multiplicity pattern $(2,1,1)$.  For a fixed repeated weight, there are
$\binom{6}{2}=15$ unordered pairs of other weights.  Three pairs lie
with the repeated weight in a two-dimensional subspace, so the remaining
$12$ pairs span $E$ with it.  This gives $7\cdot12=84$ monomials.

If all four weights are distinct, they necessarily span $E$, since a
two-dimensional subspace has only three nonzero elements.  Their unique
linear dependence either involves three weights or all four weights.
In the first case, three weights are the nonzero elements of a
two-dimensional subspace and the fourth lies outside it; there are
$7\cdot4=28$ such sets.  The remaining
$\binom{7}{4}-28=7$ sets have a four-term dependence.  Thus, up to the
$\Gamma$-action, every faithful degree-four monomial has one of the
three forms
\[
 M_1=\rho_1^2\rho_2\rho_3,\qquad
 M_2=\rho_1\rho_2\rho_{12}\rho_3,\qquad
 M_3=\rho_1\rho_2\rho_3\rho_{123}.
\]
Their orbits have sizes $84$, $28$, and $7$, respectively; these are
the types $T_1,T_2,T_3$.  Hence there are $119$ faithful degree-four
monomials in total.

By the Theorem~\ref{A} and Stong's injectivity theorem, it suffices to
compute the space of coefficient systems satisfying the evenness
conditions.  Fix one nonzero $\rho$.  The six nonzero weights
different from $\rho$ form three pairs
\[
 \{\beta,\beta+\rho\},
\]
whose members have the same reduction modulo $\rho$.  Thus the
$d=2$ classes are obtained by choosing two different pairs, while the
$d=1$ classes have quotient patterns $1+1+1$ or $2+1$.  Consequently,
for each $\rho$ there are three $d=2$ classes and seven $d=1$ classes.
For $d=2$, the four choices of actual weights give the four one-tuple
conditions in addition to the empty-tuple condition.  For $d=1$, only
the empty-tuple condition occurs.  By
Lemma~\ref{lem:tuple-test-reduction}, these are all the equations needed.
We use the following three representatives for the three parts of the
calculation; they will also be realized in the geometric constructions
that follow.  Set
\begin{align}
 \Lambda_1={}&
 \rho_1\rho_2\rho_3\rho_{123}
 +\rho_1\rho_{12}\rho_{23}\rho_3
 +\rho_2\rho_{12}\rho_{23}\rho_{123},
 \label{eq:Lambda1}\\
 \Lambda_2={}&
 \rho_1\rho_2\rho_{12}\rho_3
 +\rho_1\rho_2\rho_{12}\rho_{13}
 +\rho_1\rho_2\rho_{12}\rho_{23}
 +\rho_1\rho_2\rho_{12}\rho_{123}
 +\rho_3\rho_{13}\rho_{23}\rho_{123},
 \label{eq:Lambda2}\\
 \Lambda_3={}&
 \rho_1^2\rho_2\rho_3
 +\rho_1^2\rho_{12}\rho_3
 +\rho_1^2\rho_2\rho_{13}
 +\rho_1^2\rho_{12}\rho_{13}
 \notag\\
 &+\rho_1\rho_2\rho_{12}\rho_{23}
 +\rho_1\rho_2\rho_{12}\rho_{123}
 +\rho_1\rho_3\rho_{13}\rho_{23}
 +\rho_1\rho_3\rho_{13}\rho_{123}.
 \label{eq:Lambda3}
\end{align}

The support of $\Lambda_1$ is a Fano line among the seven monomials of
type $T_3$, so its orbit has seven elements.  The $T_2$-part of
$\Lambda_2$ is determined by one of the seven two-dimensional
subspaces of $E$, so its orbit also has seven elements.  Finally, the
$T_1$-part of $\Lambda_3$ is determined by its repeated weight and
one of the three unordered bases of the corresponding quotient.
Hence its orbit has $7\cdot3=21$ elements.

\begin{theorem}\label{thm:Z4-G3}
As a vector space over $\mathbb Z_2$,
\[
 \dim_{\mathbb Z_2}\mathcal Z_4(G_3)=32.
\]
\end{theorem}

\begin{proof}
By Theorem~\ref{A}, it remains to determine the coefficient systems
satisfying these parity equations.  We examine the coefficients of the
three monomial types in turn.

First, we verify the evenness condition for all seven nonzero
weights.  The table records the cardinalities of the $d=1$ classes:
\[
\begin{array}{c|ccccccc}
 & \rho_1&\rho_2&\rho_3&\rho_{12}&\rho_{13}&\rho_{23}&\rho_{123}\\ \hline
\Lambda_1&2&2&2&2&--&2&2\\
\Lambda_2&2+2&2+2&2&2+2&2&2&2\\
\Lambda_3&2+2&2+2&2+2&2+2&2+2&2&2
\end{array}
\]
Here $2+2$ means two classes, each of cardinality $2$, while a dash
means that no class occurs.
For $\Lambda_3$ and $\rho_1$, there is in addition one $d=2$ class
consisting of the four terms containing $\rho_1^2$.  Its empty-tuple
multiplicity is $4$, and, after removing $\rho_1^2$, each of
\[
 \rho_2,\quad \rho_{12},\quad \rho_3,\quad \rho_{13}
\]
occurs twice.  This checks every required tuple for this class as well.
Consequently, the table and the preceding $d=2$ calculation verify the
evenness condition for each of the three displayed representatives and
every nonzero $\rho$.  Since a change of basis permutes the nonzero
weights, the mod $\rho$ classes, and their tuple multiplicities, the
same verification applies to all changes of basis of
$\Lambda_1,\Lambda_2,$ and $\Lambda_3$.  Thus these representatives
and all their changes of basis satisfy the evenness condition.

We now count the free coefficients in the order $T_1,T_2,T_3$.

\emph{Repeated-weight blocks.}
For a fixed nonzero $\rho$, the six other nonzero weights form
three pairs modulo $\rho$.  Choosing two different pairs gives a
$d=2$ class consisting of four $T_1$ monomials.  For example, when
$\rho=\rho_1$, one such block is
\[
 \rho_1^2\rho_2\rho_3,\qquad
 \rho_1^2\rho_{12}\rho_3,\qquad
 \rho_1^2\rho_2\rho_{13},\qquad
 \rho_1^2\rho_{12}\rho_{13}.
\]
The other blocks are obtained from the other choices of two pairs.
For the displayed block, the $\rho_2$- and $\rho_{12}$-equations equate
the first coefficient with the third and the second with the fourth,
while the $\rho_3$-equation equates the first two.  Thus all four
coefficients are equal; the empty-tuple equation is then automatic.
The same argument applies to every such block.  As
$\rho$ and the two pairs vary, these are $7\cdot3=21$ disjoint
four-term blocks, covering all $84$ monomials of type $T_1$.  Thus the
$T_1$ coefficients contribute $21$ free parameters.  The $T_1$-parts
of the $21$ changes of basis of $\Lambda_3$ give the corresponding
leading terms.  After subtracting the corresponding admissible
systems, all coefficients of type $T_1$ may be taken to be zero.

\emph{Two-plane blocks.}
For each $\rho$, six of the seven $d=1$ classes contain four $T_1$
monomials and two $T_2$ monomials.  After the $T_1$ part has been
removed, their equations become pair-relations among the $T_2$
coefficients.  Over all $\rho$, these $42$ relations split the
$28$ type $T_2$ monomials into seven connected components of four
terms.  The components are indexed by the seven two-dimensional
subspaces, with the four choices of the weight outside the subspace.
For example, the block associated with
$\{\rho_1,\rho_2,\rho_{12}\}$ is
\[
 \rho_1\rho_2\rho_{12}\rho_3,\qquad
 \rho_1\rho_2\rho_{12}\rho_{13},\qquad
 \rho_1\rho_2\rho_{12}\rho_{23},\qquad
 \rho_1\rho_2\rho_{12}\rho_{123}.
\]
For this block, the equations modulo $\rho_1$ equate the first
coefficient with the second and the third with the fourth.  The
equations modulo $\rho_2$ equate the first coefficient with the third
and the second with the fourth.  Hence all four coefficients are equal.
Every other two-dimensional subspace is obtained by a change of basis,
so the same conclusion holds for all seven blocks.
Hence the $T_2$ coefficients contribute seven free parameters.  The
corresponding seven changes of basis of $\Lambda_2$ satisfy all the
evenness equations; in each case, the $T_2$ block together with the
remaining $T_3$ term satisfies the mixed $d=1$ equations.
After subtracting them, we may likewise take all coefficients of type
$T_2$ to be zero.

\goodbreak
\emph{Fano-plane equations.}
The remaining seven equations now involve only the $T_3$ coefficients.
They are the incidence equations for the complements of the seven
Fano lines.  Under this indexing, the equations corresponding to
$\rho_1,\rho_2,\rho_3$ are independent, while those corresponding to
$\rho_{12},\rho_{13},\rho_{23},\rho_{123}$ are respectively their
pairwise sums and their total sum.  Thus these equations have rank $3$
and leave a
$4$-dimensional solution space.  This is the span of the $T_3$-parts
of the changes of basis of $\Lambda_1$: these are the seven Fano-line
incidence vectors, which satisfy the equations and have rank $4$.

The three steps show that every coefficient system satisfying the
evenness condition is a linear combination of $21+7+4=32$ such
systems.  They are independent: in any linear relation, the disjoint
$T_1$-supports first force the coefficients of the $21$ changes of
basis of $\Lambda_3$ to vanish.  The disjoint $T_2$-supports then force
the seven coefficients belonging to $\Lambda_2$ to vanish, and the
remaining $T_3$-space has dimension $4$.  Hence the $32$ systems form a
basis of the space of admissible coefficient systems.  Therefore, by
the Theorem~\ref{A}
and Stong's injectivity theorem,
\[
 \dim_{\mathbb Z_2}\mathcal Z_4(G_3)=32.
\]
\end{proof}

\subsection{Geometric representatives}

The preceding calculation separates the generators into three orbit
types.  We now realize the corresponding data
$\Lambda_1,\Lambda_2,\Lambda_3$ of
\eqref{eq:Lambda1}--\eqref{eq:Lambda3} by explicit equivariant
manifolds.  The leading parts of $\Lambda_1$, $\Lambda_2$, and
$\Lambda_3$ belong respectively to the square-free, three-in-a-plane,
and repeated-weight types.  Their images under changes of basis give
the $4$, $7$, and $21$ directions found above.  We begin with the only
one involving a repeated weight.

\subsubsection*{\texorpdfstring{$\bullet$}{\textbullet} The generator $\Lambda_3$}
Let $\gamma$ be the canonical line bundle over $\mathbb RP^1$ and set
\[
X=\mathbb RP(\gamma \oplus \gamma \oplus \gamma  \oplus\underline{\mathbb R}).
\]
This is the small cover over $\Delta^3\times\Delta^1$ \cite{DJ} whose six
facet vectors are
\[
e_1,\quad e_2,\quad e_3,\quad e_1+e_2+e_3,\quad
e_4,\quad e_1+e_2+e_3+e_4.
\]
Equivalently, its characteristic matrix is
\[
\begin{pmatrix}
1&0&0&1&0&1\\
0&1&0&1&0&1\\
0&0&1&1&0&1\\
0&0&0&0&1&1
\end{pmatrix}.
\]
Let $\beta_1,\ldots,\beta_4$ be the standard basis of
$\Hom(G_4,\mathbb Z_2)$, and
restrict the $G_4$-action to
\[
H=\ker(\beta_1+\beta_2+\beta_4)\cong G_3,\qquad
\rho_i=\beta_i|_H\quad(1\leq i\leq3).
\]
Then $\beta_4|_H=\rho_{12}$.  The eight fixed points have tangent
representations
\[
\begin{gathered}
\rho_{12}^2\rho_1\rho_{13},\quad
\rho_{12}^2\rho_2\rho_{13},\quad
\rho_{12}^2\rho_2\rho_{23},\quad
\rho_{12}^2\rho_1\rho_{23},\\
\rho_{12}\rho_{13}\rho_{23}\rho_3,\quad
\rho_{12}\rho_{13}\rho_{23}\rho_{123},\quad
\rho_1\rho_2\rho_{12}\rho_3,\quad
\rho_1\rho_2\rho_{12}\rho_{123}.
\end{gathered}
\]
No restricted tangent weight is trivial, so these are all the
$H$-fixed points: a positive-dimensional $H$-fixed stratum of a small
cover would contain a vertex and give a trivial restricted tangent
weight there.  Applying the automorphism
$\rho_1\mapsto\rho_2$, $\rho_2\mapsto\rho_{12}$, and
$\rho_3\mapsto\rho_{23}$ to the displayed fixed point polynomial
gives $\Lambda_3$.

The following labeled graph records the corresponding eight tangent
representations.
\begin{center}
\begin{tikzpicture}[scale=.9]
  \draw[thick] (-1,0.5) rectangle (1,-0.5);

  \draw[thick] (1,0.5)..controls(2,0.85)..(2.5,1);
  \draw[thick] (-1,0.5)..controls(-2,0.85)..(-2.5,1);
  \draw[thick] (1,-0.5)..controls(2,-0.85)..(2.5,-1);
  \draw[thick] (-1,-0.5)..controls(-2,-0.85)..(-2.5,-1);

  \draw[thick] (1,0.5)--(2.5,-1);
  \draw[thick] (-1,0.5)--(-2.5,-1);
  \draw[thick] (1,-0.5)--(2.5,1);
  \draw[thick] (-1,-0.5)--(-2.5,1);

  \draw[thick] (2.5,1)..controls(2.6,0)..(2.5,-1);
  \draw[thick] (-2.5,1)..controls(-2.6,0)..(-2.5,-1);

  \draw[thick] (2.5,1)..controls(0,1.3)..(-2.5,1);
  \draw[thick] (2.5,-1)..controls(0,-1.3)..(-2.5,-1);

  \filldraw (1,0.5) circle (.1) (-1,0.5) circle (.1)
  (1,-0.5) circle (.1) (-1,-0.5) circle (.1)
  (2.5,1) circle (.1) (-2.5,1) circle (.1)
  (2.5,-1) circle (.1) (-2.5,-1) circle (.1);

  \node[above] at (0,0.4) {$\rho_1$};
  \node[below] at (0,-0.4) {$\rho_1$};
  \node[above] at (0,1.2) {$\rho_1$};
  \node[below] at (0,-1.2) {$\rho_1$};

  \node[left] at (-0.4,0) {$\rho_1$};
  \node[right] at (0.5,0) {$\rho_1$};
  \node[left] at (-2.6,0) {$\rho_3$};
  \node[right] at (2.6,0) {$\rho_{13}$};

  \node[above] at (-1.5,0.6) {$\rho_2$};
  \node[above] at (1.5,0.6) {$\rho_2$};
  \node[below] at (-1.5,-0.65) {$\rho_{123}$};
  \node[below] at (1.5,-0.63) {$\rho_{23}$};

  \node[below] at (-2.2,0.63) {$\rho_{12}$};
  \node[below] at (2.2,0.63) {$\rho_{12}$};
  \node[above] at (-2.2,-0.6) {$\rho_{23}$};
  \node[above] at (2.23,-0.56) {$\rho_{123}$};
\end{tikzpicture}
\end{center}

\subsubsection*{\texorpdfstring{$\bullet$}{\textbullet} The generator $\Lambda_2$}
Let $G_4$ act linearly on $\mathbb RP^4$ by changing the signs of the
last four homogeneous coordinates:
\[
[x_0:x_1:x_2:x_3:x_4]\longmapsto
[x_0:(-1)^{g_1}x_1:(-1)^{g_2}x_2:
(-1)^{g_3}x_3:(-1)^{g_4}x_4].
\]
Restrict to the same subgroup $H$.  The coordinate weights
$0,\beta_1,\ldots,\beta_4$ have pairwise distinct restrictions to
$H$, so the five coordinate points are the whole fixed point set.
At a coordinate fixed point, the tangent representation is obtained
from the four remaining coordinate directions.  Hence the five tangent
representations are
\begin{align*}
&\rho_1\rho_2\rho_{12}\rho_3,\quad
\rho_1\rho_2\rho_{12}\rho_{13},\quad
\rho_1\rho_2\rho_{12}\rho_{23},\\
&\rho_1\rho_2\rho_{12}\rho_{123},\quad
\rho_3\rho_{13}\rho_{23}\rho_{123}.
\end{align*}
Their sum is exactly $\Lambda_2$.

The associated labeled graph is shown below.
\begin{center}
\begin{tikzpicture}[scale=.9]
  \draw[thick] (-1,0)--(1,0)--(1,-2)--(-1,-2)--cycle;
  \draw[thick] (-1,0)--(1,-2);
  \draw[thick] (1,0)--(-1,-2);
  \draw[thick] (-1,0)--(0,1);
  \draw[thick] (1,0)--(0,1);
  \draw[thick] (-1,-2)..controls (-2.5, -0.3) and (-2,0.5) .. (0,1);
  \draw[thick] (1,-2)..controls (2.5, -0.3) and (2,0.5) .. (0,1);
  \filldraw (-1,0) circle (.1) (1,0) circle (.1)
  (1,-2) circle (.1) (-1,-2) circle (.1) (0,1) circle (.1);

  \node[right] at (1.8,-0.5) {$\rho_{123}$};
  \node[left] at (-1.8,-0.5) {$\rho_3$};
  \node[left] at (-0.5,0.45) {$\rho_{13}$};
  \node[right] at (0.5,0.5) {$\rho_{23}$};

  \node[left] at (-0.9,-1) {$\rho_1$};
  \node[right] at (0.9,-1) {$\rho_1$};

  \node[above] at (0,-0.1) {$\rho_2$};
  \node[below] at (0,-1.9) {$\rho_2$};

  \node[left] at (-0.2,-0.78) {$\rho_{12}$};
  \node[right] at (0.2,-0.78) {$\rho_{12}$};
\end{tikzpicture}
\end{center}

\subsubsection*{\texorpdfstring{$\bullet$}{\textbullet} The generator $\Lambda_1$}
The standard $G_3$-twisted involution over
$\mathbb RP^2\times\mathbb RP^2$ is represented as follows:
\[
 a_{g_1,g_2}[u_0,u_1,u_2]
 =[u_0,(-1)^{g_1}u_1,(-1)^{g_2}u_2],
\]
and
\[
 (g_1,g_2,g_3)\cdot([z],[w])=
 \begin{cases}
  \bigl(a_{g_1,g_2}[z],a_{g_1,g_2}[w]\bigr),&g_3=0,\\
  \bigl(a_{g_1,g_2}[w],a_{g_1,g_2}[z]\bigr),&g_3=1.
 \end{cases}
\]

The fixed points are the three diagonal pairs of coordinate points,
and their fixed point data are
\[
\rho_1\rho_2\rho_{13}\rho_{23},\qquad
\rho_1\rho_{12}\rho_{13}\rho_{123},\qquad
\rho_2\rho_{12}\rho_{23}\rho_{123}.
\]
The automorphism fixing $\rho_1,\rho_2$ and sending
$\rho_3$ to $\rho_{13}$ transforms their sum into $\Lambda_1$.

This generator is represented by the following labeled graph.
\begin{center}
\begin{tikzpicture}[scale=.9]
  \draw[thick] (-1.6,0)..controls (0,0.4).. (1.6,0);
  \draw[thick] (-1.6,0)..controls (0,-0.4).. (1.6,0);

  \draw[thick] (-1.6,0)..controls (-1,1.5).. (0,2);
  \draw[thick] (-1.6,0)..controls (-1,0.3).. (0,2);

  \draw[thick] (1.6,0)..controls (1,1.5).. (0,2);
  \draw[thick] (1.6,0)..controls (1,0.3).. (0,2);

  \filldraw (-1.6,0) circle (.1) (1.6,0) circle (.1) (0,2) circle (.1);

  \node[right] at (-0.7,0.9) {$\rho_1$};
  \node[left] at (-1.2,1) {$\rho_{13}$};
  \node[left] at (0.7,0.9) {$\rho_2$};
  \node[right] at (1.2,1) {$\rho_{23}$};

  \node[above] at (0,0.2) {$\rho_{12}$};
  \node[below] at (0,-0.2) {$\rho_{123}$};
\end{tikzpicture}
\end{center}

\subsection{A basis of $\Zc_4(G_3)$}\label{subsec:complete-list}

The $\Gamma$-orbits of $\Lambda_1$, $\Lambda_2$, and $\Lambda_3$ are
obtained from the three actions above by precomposition with
automorphisms of $G_3$.  They have sizes $7$, $7$, and $21$, while
their spans have dimensions $4$, $7$, and $21$, respectively, and
together generate $\mathcal Z_4(G_3)$.  Thus a basis is obtained by
choosing four independent elements from the first orbit and every
element from the other two.  Its $32$ elements are represented by four
actions on $\mathbb RP^2\times\mathbb RP^2$, seven on $\mathbb RP^4$,
and twenty-one on
$\mathbb RP(\gamma\oplus\gamma\oplus\gamma\oplus
\underline{\mathbb R})$.  Their fixed-point data are listed below, and
their linear independence is verified in the next subsection.  Here
$\rho_{ij}=\rho_i+\rho_j$ and $\rho_{123}=\rho_1+\rho_2+\rho_3$.

\medskip
\noindent\textbf{Square-free type (4 generators):}
\[
\begin{aligned}
f_{1,1} &= \rho_1\rho_2\rho_3\rho_{123}
        + \rho_1\rho_{12}\rho_{23}\rho_3
        + \rho_2\rho_{12}\rho_{23}\rho_{123},\\
f_{1,2} &= \rho_1\rho_2\rho_3\rho_{123}
        + \rho_1\rho_2\rho_{13}\rho_{23}
        + \rho_3\rho_{13}\rho_{23}\rho_{123},\\
f_{1,3} &= \rho_1\rho_2\rho_{13}\rho_{23}
        + \rho_1\rho_{12}\rho_{23}\rho_3
        + \rho_2\rho_{12}\rho_{13}\rho_3,\\
f_{1,4} &= \rho_1\rho_2\rho_{13}\rho_{23}
        + \rho_1\rho_{12}\rho_{13}\rho_{123}
        + \rho_2\rho_{12}\rho_{23}\rho_{123}.
\end{aligned}
\]

\noindent\textbf{Three-in-a-plane type (7 generators):}
\[
\begin{aligned}
f_{2,1} &= \rho_1\rho_2\rho_{12}\rho_3 + \rho_1\rho_2\rho_{12}\rho_{13}
        + \rho_1\rho_2\rho_{12}\rho_{23} + \rho_1\rho_2\rho_{12}\rho_{123}
        + \rho_{13}\rho_{23}\rho_3\rho_{123},\\
f_{2,2} &= \rho_1\rho_2\rho_3\rho_{13} + \rho_1\rho_{12}\rho_{13}\rho_3
        + \rho_1\rho_{13}\rho_{23}\rho_3 + \rho_1\rho_{13}\rho_3\rho_{123}
        + \rho_{12}\rho_2\rho_{23}\rho_{123},\\
f_{2,3} &= \rho_1\rho_2\rho_3\rho_{23} + \rho_2\rho_3\rho_{23}\rho_{12}
        + \rho_2\rho_3\rho_{23}\rho_{13} + \rho_2\rho_3\rho_{23}\rho_{123}
        + \rho_1\rho_{12}\rho_{13}\rho_{123},\\
f_{2,4} &= \rho_1\rho_{12}\rho_{13}\rho_{23} + \rho_{12}\rho_{13}\rho_{23}\rho_2
        + \rho_{12}\rho_{13}\rho_{23}\rho_3
        + \rho_{12}\rho_{13}\rho_{23}\rho_{123}
        + \rho_1\rho_2\rho_3\rho_{123},\\
f_{2,5} &= \rho_{123}\rho_1\rho_{23}\rho_2 + \rho_{123}\rho_1\rho_{23}\rho_3
        + \rho_{123}\rho_1\rho_{23}\rho_{12}
        + \rho_{123}\rho_1\rho_{23}\rho_{13}
        + \rho_2\rho_3\rho_{12}\rho_{13},\\
f_{2,6} &= \rho_{123}\rho_2\rho_{13}\rho_1 + \rho_{123}\rho_2\rho_{13}\rho_3
        + \rho_{123}\rho_2\rho_{13}\rho_{12}
        + \rho_{123}\rho_2\rho_{13}\rho_{23}
        + \rho_1\rho_3\rho_{12}\rho_{23},\\
f_{2,7} &= \rho_{123}\rho_3\rho_{12}\rho_1 + \rho_{123}\rho_3\rho_{12}\rho_2
        + \rho_{123}\rho_3\rho_{12}\rho_{23}
        + \rho_{123}\rho_3\rho_{12}\rho_{13}
        + \rho_1\rho_2\rho_{13}\rho_{23}.
\end{aligned}
\]

\noindent\textbf{Repeated-weight type (21 generators):}
\[
\resizebox{\textwidth}{!}{\small$\begin{aligned}
&\rho_1^2\rho_2\rho_3+ \rho_1^2\rho_{12}\rho_3+\rho_1^2\rho_2\rho_{13}
 + \rho_1^2\rho_{12}\rho_{13}
 + \rho_1\rho_2\rho_{12}\rho_{23}+\rho_1\rho_2\rho_{12}\rho_{123}
 + \rho_1\rho_3\rho_{13}\rho_{23}+\rho_1\rho_3\rho_{13}\rho_{123},\\
&\rho_1^2\rho_2\rho_{23}+ \rho_1^2\rho_2\rho_{123}
 + \rho_1^2\rho_{12}\rho_{23}+ \rho_1^2\rho_{12}\rho_{123}
 + \rho_1\rho_2\rho_{12}\rho_3+ \rho_1\rho_2\rho_{12}\rho_{13}
 + \rho_1\rho_{23}\rho_{123}\rho_3+ \rho_1\rho_{23}\rho_{123}\rho_{13},\\
&\rho_1^2\rho_{23}\rho_3+ \rho_1^2\rho_{23}\rho_{13}
 + \rho_1^2\rho_{123}\rho_3+ \rho_1^2\rho_{123}\rho_{13}
 + \rho_1\rho_{23}\rho_{123}\rho_2+ \rho_1\rho_{23}\rho_{123}\rho_{12}
 + \rho_1\rho_3\rho_{13}\rho_2+\rho_1\rho_3\rho_{13}\rho_{12},\\
&\rho_2^2\rho_1\rho_3+ \rho_2^2\rho_1\rho_{23}
 + \rho_2^2\rho_{12}\rho_3+ \rho_2^2\rho_{12}\rho_{23}
 + \rho_1\rho_2\rho_{12}\rho_{13}+ \rho_1\rho_2\rho_{12}\rho_{123}
 + \rho_2\rho_3\rho_{23}\rho_{13}+ \rho_2\rho_3\rho_{23}\rho_{123},\\
&\rho_2^2\rho_1\rho_{13}+ \rho_2^2\rho_1\rho_{123}
 + \rho_2^2\rho_{12}\rho_{13}+ \rho_2^2\rho_{12}\rho_{123}
 + \rho_1\rho_2\rho_{12}\rho_3+ \rho_1\rho_2\rho_{12}\rho_{23}
 + \rho_2\rho_{13}\rho_{123}\rho_3+ \rho_2\rho_{13}\rho_{123}\rho_{23},\\
&\rho_2^2\rho_{13}\rho_3+ \rho_2^2\rho_{13}\rho_{23}
 + \rho_2^2\rho_{123}\rho_3+ \rho_2^2\rho_{123}\rho_{23}
 + \rho_{13}\rho_2\rho_{123}\rho_1+ \rho_{13}\rho_2\rho_{123}\rho_{12}
 + \rho_2\rho_3\rho_{23}\rho_1+ \rho_2\rho_3\rho_{23}\rho_{12},\\
&\rho_3^2\rho_2\rho_1+ \rho_3^2\rho_2\rho_{13}
 + \rho_3^2\rho_{23}\rho_1+ \rho_3^2\rho_{23}\rho_{13}
 + \rho_3\rho_2\rho_{23}\rho_{12}+ \rho_3\rho_2\rho_{23}\rho_{123}
 + \rho_3\rho_1\rho_{13}\rho_{12}+ \rho_3\rho_1\rho_{13}\rho_{123},\\
&\rho_3^2\rho_2\rho_{12}+ \rho_3^2\rho_2\rho_{123}
 + \rho_3^2\rho_{23}\rho_{12}+ \rho_3^2\rho_{23}\rho_{123}
 + \rho_3\rho_2\rho_{23}\rho_1+ \rho_3\rho_2\rho_{23}\rho_{13}
 + \rho_3\rho_{12}\rho_{123}\rho_1+ \rho_3\rho_{12}\rho_{123}\rho_{13},\\
&\rho_3^2\rho_{12}\rho_1+ \rho_3^2\rho_{12}\rho_{13}
 + \rho_3^2\rho_{123}\rho_1+ \rho_3^2\rho_{123}\rho_{13}
 + \rho_3\rho_{12}\rho_{123}\rho_2+ \rho_3\rho_{12}\rho_{123}\rho_{23}
 + \rho_3\rho_1\rho_{13}\rho_2+ \rho_3\rho_1\rho_{13}\rho_{23},\\
&\rho_{12}^2\rho_1\rho_3+ \rho_{12}^2\rho_1\rho_{123}
 + \rho_{12}^2\rho_2\rho_3+ \rho_{12}^2\rho_2\rho_{123}
 + \rho_1\rho_2\rho_{12}\rho_{13}+ \rho_1\rho_2\rho_{12}\rho_{23}
 + \rho_{12}\rho_3\rho_{123}\rho_{13}+ \rho_{12}\rho_3\rho_{123}\rho_{23},\\
&\rho_{12}^2\rho_1\rho_{13}+ \rho_{12}^2\rho_1\rho_{23}
 + \rho_{12}^2\rho_2\rho_{13}+ \rho_{12}^2\rho_2\rho_{23}
 + \rho_1\rho_2\rho_{12}\rho_3+ \rho_1\rho_2\rho_{12}\rho_{123}
 + \rho_{12}\rho_{13}\rho_{23}\rho_3+ \rho_{12}\rho_{13}\rho_{23}\rho_{123},\\
&\rho_{12}^2\rho_{13}\rho_3+ \rho_{12}^2\rho_{13}\rho_{123}
 + \rho_{12}^2\rho_{23}\rho_3+ \rho_{12}^2\rho_{23}\rho_{123}
 + \rho_{13}\rho_{23}\rho_{12}\rho_1+ \rho_{13}\rho_{23}\rho_{12}\rho_2
 + \rho_{12}\rho_3\rho_{123}\rho_1+ \rho_{12}\rho_3\rho_{123}\rho_2,\\
&\rho_{13}^2\rho_2\rho_3+ \rho_{13}^2\rho_2\rho_1
 + \rho_{13}^2\rho_{123}\rho_3+ \rho_{13}^2\rho_{123}\rho_1
 + \rho_{13}\rho_2\rho_{123}\rho_{23}+ \rho_{13}\rho_2\rho_{123}\rho_{12}
 + \rho_{13}\rho_3\rho_1\rho_{23}+ \rho_{13}\rho_3\rho_1\rho_{12},\\
&\rho_{13}^2\rho_2\rho_{23}+ \rho_{13}^2\rho_2\rho_{12}
 + \rho_{13}^2\rho_{123}\rho_{23}+ \rho_{13}^2\rho_{123}\rho_{12}
 + \rho_{13}\rho_2\rho_{123}\rho_3+ \rho_{13}\rho_2\rho_{123}\rho_1
 + \rho_{13}\rho_{23}\rho_{12}\rho_3+\rho_{13}\rho_{23}\rho_{12}\rho_1,\\
&\rho_{13}^2\rho_{23}\rho_3+ \rho_{13}^2\rho_{23}\rho_1
 + \rho_{13}^2\rho_{12}\rho_3+ \rho_{13}^2\rho_{12}\rho_1
 + \rho_{13}\rho_{23}\rho_{12}\rho_2+ \rho_{13}\rho_{23}\rho_{12}\rho_{123}
 + \rho_{13}\rho_3\rho_1\rho_2+ \rho_{13}\rho_3\rho_1\rho_{123},\\
&\rho_{23}^2\rho_1\rho_3+ \rho_{23}^2\rho_1\rho_2
 + \rho_{23}^2\rho_{123}\rho_3+ \rho_{23}^2\rho_{123}\rho_2
 + \rho_1\rho_{23}\rho_{123}\rho_{13}+ \rho_1\rho_{23}\rho_{123}\rho_{12}
 + \rho_{23}\rho_3\rho_2\rho_{13}+ \rho_{23}\rho_3\rho_2\rho_{12},\\
&\rho_{23}^2\rho_1\rho_{13}+ \rho_{23}^2\rho_1\rho_{12}
 + \rho_{23}^2\rho_{123}\rho_{13}+ \rho_{23}^2\rho_{123}\rho_{12}
 + \rho_1\rho_{23}\rho_{123}\rho_3+ \rho_1\rho_{23}\rho_{123}\rho_2
 + \rho_{23}\rho_{13}\rho_{12}\rho_3+ \rho_{23}\rho_{13}\rho_{12}\rho_2,\\
&\rho_{23}^2\rho_{13}\rho_3+ \rho_{23}^2\rho_{13}\rho_2
 + \rho_{23}^2\rho_{12}\rho_3+ \rho_{23}^2\rho_{12}\rho_2
 + \rho_{13}\rho_{23}\rho_{12}\rho_1+ \rho_{13}\rho_{23}\rho_{12}\rho_{123}
 + \rho_{23}\rho_3\rho_2\rho_1+ \rho_{23}\rho_3\rho_2\rho_{123},\\
&\rho_{123}^2\rho_2\rho_3+ \rho_{123}^2\rho_2\rho_{12}
 + \rho_{123}^2\rho_{13}\rho_3+ \rho_{123}^2\rho_{13}\rho_{12}
 + \rho_{123}\rho_2\rho_{13}\rho_{23}+ \rho_{123}\rho_2\rho_{13}\rho_1
 + \rho_{123}\rho_3\rho_{12}\rho_{23}+ \rho_{123}\rho_3\rho_{12}\rho_1,\\
&\rho_{123}^2\rho_2\rho_{23}+ \rho_{123}^2\rho_2\rho_1
 + \rho_{123}^2\rho_{13}\rho_{23}+ \rho_{123}^2\rho_{13}\rho_1
 + \rho_{123}\rho_2\rho_{13}\rho_3+ \rho_{123}\rho_2\rho_{13}\rho_{12}
 + \rho_{123}\rho_{23}\rho_1\rho_3+ \rho_{123}\rho_{23}\rho_1\rho_{12},\\
&\rho_{123}^2\rho_{23}\rho_3+ \rho_{123}^2\rho_{23}\rho_{12}
 + \rho_{123}^2\rho_1\rho_3+ \rho_{123}^2\rho_1\rho_{12}
 + \rho_{123}\rho_{23}\rho_1\rho_2+ \rho_{123}\rho_{23}\rho_1\rho_{13}
 + \rho_{123}\rho_3\rho_{12}\rho_2+\rho_{123}\rho_3\rho_{12}\rho_{13}.
\end{aligned}$}
\]


The $32$ fixed point data listed above are linearly independent over
$\mathbb Z_2$.  Indeed, the repeated-weight parts of the last $21$
data have pairwise disjoint supports, so their coefficients in any
linear relation must vanish.  The three-in-a-plane parts of the next
seven data also have pairwise disjoint supports, so their coefficients
must then vanish.  Finally, the first four square-free data are four
independent Fano-line incidence vectors.  Thus all $32$ coefficients
vanish.  Together with
$\dim_{\mathbb Z_2}\mathcal Z_4(G_3)=32$ from
Theorem~\ref{thm:Z4-G3}, these data form a basis of
$\mathcal Z_4(G_3)$.

\bibliographystyle{abbrv}

\end{document}